\documentclass[12pt]{article}
\usepackage{mathrsfs}
\usepackage{amssymb,a4}
\usepackage{amsmath,amsfonts,amssymb,amsthm}
\newcommand{\al}{\alpha}

\def\Res{{\rm Res}}
\def\wt{{\rm wt}}

\def\be{\beta}

\newcommand{\la}{\lambda}

\def\C{{\mathbb C}}

\def\Z{{\mathbb Z}}

\def\N{{\mathbb N}}
\def\1{{\bf 1}}

\def \wt{{\rm wt}}

\def \Res{{\rm Res}}

\def \mod{{\rm mod}}

\def \<{\langle}
\def \>{\rangle}

\def \pf{\noindent {\bf Proof: \,}}
\def\theequation{5.\arabic{equation}}

\renewcommand{\theequation}{\thesection.\arabic{equation}}

\newtheorem{theorem}{Theorem}[section]
\newtheorem{prop}[theorem]{Proposition}
\newtheorem{lem}[theorem]{Lemma}

\newtheorem{remark}[theorem]{Remark}
\theoremstyle{definition}

\begin{document}
\begin{center}
{\Large {\bf A characterization of the rational vertex operator
algebra $V_{\Z\al}^{+}$}: II} \\

\vspace{0.5cm} Chongying Dong\footnote{Supported by NSF grants and a
Faculty research grant from  the University of California at Santa
Cruz.}
\\
Department of Mathematics,  University of California, Santa Cruz, CA 95064 \\
\vspace{.1 cm} Cuipo Jiang\footnote{Supported  by China NSF grants
10931006,10871125, the RFDP of China(20100073110052), and the
Innovation Program of Shanghai Municipal
Education Commission (11ZZ18).}\\
 Department of Mathematics, Shanghai Jiaotong University, Shanghai 200240 China
\end{center}
\hspace{1cm}
\begin{abstract} A characterization of vertex operator algebra $V_L^+$ for any rank one positive definite
even lattice $L$ is given in terms of dimensions of homogeneous subspaces of small weights. This result reduces the
classification of rational vertex operator algebras of central charge $1$ to the characterization of three
 vertex operator algebras in the $E$-series of central charge one.

2000MSC:17B69
\end{abstract}

\section{Introduction}

It is a well known conjecture in the theory of vertex operator
algebra that any rational vertex operator algebras with central
charge $c=1$ is isomorphic to $V_L,$ or $V_L^+$ or $V_{\Z\alpha}^G$
where $L$ is a rank one positive definite even lattice,
$(\alpha,\alpha)=2$ and $G=A_4,S_4,A_5$ is a subgroup of $SO(3)$ in
the $E$-series (cf. \cite{K}). The characterization of $V_L$ for any
positive definite even lattice was established in \cite{DM1} in
terms the rank of $V_1,$ central charge and effective central
charge. A characterization of $V_{\Z\beta}^+$ for $(\beta,\beta)=4$
 in terms of $\dim V_{2}$ is given in \cite{ZD}, \cite{DJ1}.
 Characterization for all $V_{\Z\beta}^+$ with $(\beta,\beta)/2$ not
being a perfect square  in terms of dimensions of $V_i$ for $i\leq
4$ is obtained in \cite{DJ2}. In this paper we characterize
$V_{\Z\beta}$ with $(\beta,\beta)/2$ being a perfect square by
dimensions of $V_i$ for $i\leq 4.$ It remains the characterization
of $V_{\Z\alpha}^G$ for $G=A_4,S_4$ and $A_5$ for completing the
classification of rational vertex operator algebras with $c=1.$

There are two major differences between $V_{\Z\beta}^+$ and
$V_{\Z\alpha}^G.$ The first one is that
 $V_{\Z\beta}^+$ is the fixed points of rational vertex operator algebra $V_{\Z\beta}$
 under an order two automorphism and $V_{\Z\alpha}^G$ is the fixed points of rational vertex operator algebra $V_{\Z\alpha}$
  under a nonabelian group. The rationality and classification of irreducible modules of $V_{\Z\alpha}^G$ have not been achieved although
   the automorphism groups of $V_{\Z\alpha}^G$ are known
 \cite{DG}, \cite{DGR}. But this difference is not our concern in this paper.
  The second difference comes from the dimensions of weight 4 subspaces: $\dim(V_{\Z\beta}^+)_4\geq 3$ and $\dim(V_{\Z\alpha}^G)_4=2.$
 This difference inspires us  to characterize  $V_{\Z\beta}^+$ in terms of dimensions of $V_{i}$ for $i\leq 4$ in \cite{DJ2} and this paper.
 So one natural assumption for the  vertex operator algebra $V$ with $c=1$  discussed in this paper is  $\dim V_4\geq 3.$

Although the rank of $L=\Z\beta$ is one, the vertex operator algebra
$V_L^+$ is still a hard object to study. As the weight one subspace
is zero and weight two subspace is spanned by the Virasoro vector,
one can hardly use any results from the Lie algebra or Griess
algebra to obtain useful information. On the other hand, the
structure and representation theory of vertex operator algebras
$V_L^+$ and its subalgebra $M(1)^+$ have been studied extensively in
\cite{A1},  \cite{A2}, \cite{AD}, \cite{ADL}, \cite{DN1}, \cite{DN2}
and \cite{DN3}. As in \cite{DJ2} we use results from these
papers to find  a vertex operator subalgebra of $V$ isomorphic to
$M(1)^+.$ But the situation is much more complicated as $(\be,\be)/2$ is a  perfect square. We
need different ideas and methods. It is well known that
 $M(1)^+$ in the rank one case is generated by the Virasoro vector
$\omega$ and a highest weight vector $J$ of weight $4.$ The main
property of $J$ is the following: $J_3J=x+a J$ for some $x\in
L(1,0)$ which is the vertex operator subalgebra generated by
$\omega$ and some nonzero $a\in\C.$ It turns out that searching for
such $J$ in an abstract vertex operator algebra satisfying certain
assumptions is a very difficult task and involves delicate use of
the fusion rules for the vertex operator algebra $L(1,0)$ \cite{M},
\cite{DJ1}.

 We should point out that in the characterization of the lattice
vertex operator algebras $V_L$ for a positive definite even lattice
$L$ we need to use extra assumptions on the $C_2$-cofiniteness and
the effective central charge (see \cite{DLM1} for the reason).  But
we do  not need the $C_2$-cofiniteness and the effective cental
charge being one in the characterization of $V_{\Z\beta}^+$ with
$(\beta,\beta)/2>2$ not being a perfect square \cite{DJ2}.  The
situation for $(\beta,\beta)/2$ being a perfect square is totally
different. Although the effective central charge never plays any
role in this paper, the $C_2$-cofiniteness does. During the search
for vector $J$ in $V$ we cannot avoid to use the modular invariance
result from \cite{Z} where the $C_2$-cofiniteness is assumed. This
is not surprising as defining effective central charges requires the
$C_2$-cofiniteness \cite{DM1} and  the conjecture on rational vertex
operator algebras with central charge 1 is not true without assuming
the effective central charge is also one (see \cite{ZD}).

We refer the readers to \cite{AP} and \cite{X} for the related work.

This paper is organized as follows. We review the fusion rules for
the vertex operator algebra $L(1,0)$ from \cite{M} and \cite{DJ1} in
Section $2.$ We also present some results concerning Zhu algebra
\cite{Z}  and calculations on $J$ in $M(1)^+$. In Sections 3 and 4
we search for a highest weight vector $J'$ of weight 4 in an
abstract vertex operator $V$  such that $J_3'J'=x+aJ'$ for some
$x\in L(1,0)$ and a nonzero $a\in\C.$ If the space $A_4$ of highest
weight vectors of weight $4$ is one dimensional, it is trivial to
find such $J'.$ If $\dim A_4\geq 2$ this is highly nontrivial. It is
proved first that such $J'$ exists if $\dim A_4=2.$ Then it is
shown that if $\dim A_4\geq 2$ then $\dim A_4=2.$ The fusion rules
of $L(1,0)$ is used heavily here. The modular invariance of trace
functions also plays a role in this part. Section 5 is devoted to
the proof that the vertex operator subalgebra $U$ of $V$ generated
by $\omega$ and $J'$ is isomorphic to $M(1)^+$.  Section 6 gives the
main theorem: $V$ is isomorphic to $V_{\Z\beta}^+$ such that
$(\beta,\beta)=2k^2$ for some $k>1.$  A major step in this section
is to show that $V$ is a completely reducible $M(1)^+$-module with
the help of fusion rules for $M(1)^+$ and $L(1,0).$

\section{Preliminaries}
\def\theequation{2.\arabic{equation}}
\setcounter{equation}{0}

In this section we recall the fusion rules for the Virasoro vertex
operator algebra $L(1,0)$ from \cite{M} and \cite{DJ1} and for the
vertex operator algebra $M(1)^+$ of central charge $1$ from
\cite{A1}. We also discuss various results on the generator $J$ of
$M(1)^+$ following \cite{DN1}.

Let $L(c,h)$ be the highest weight irreducible module for the
Virasoro algebra with central charge $c$ and highest weight $h.$
Then $L(c,0)$ is a vertex operator algebra and each $L(c,h)$ is an
irreducible module for $L(c,0).$ In this paper we are mainly
concerned with $L(1,0)$ and its irreducible modules. First, we have
from \cite{M} and \cite{DJ1} (also see \cite{RT}):
\begin{theorem}\label{2t1} We have
\begin{equation}\label{2e1}
\dim I_{L(1,0)} \left(\hspace{-3 pt}\begin{array}{c} L(1,k^{2})\\
L(1, m^{2})\,L(1, n^{2})\end{array}\hspace{-3 pt}\right)=1,\ \
k\in{\mathbb Z}_{+},  \ |n-m|\leq k\leq n+m,\end{equation}
\begin{equation}\label{2e2}
\dim I_{L(1,0)} \left(\hspace{-3 pt}\begin{array}{c} L(1,k^{2})\\
L(1, m^{2})\,L(1, n^{2})\end{array}\hspace{-3 pt}\right)=0,\ \
k\in{\mathbb Z}_{+},  \ k<|n-m| \ {\rm or} \  k>n+m, \end{equation}
where $n,m\in{\mathbb Z}_{+}$. For $n\in{\mathbb Z}_{+}$ such that
$n\neq p^{2}$ for all $p\in{\mathbb Z}_{+}$, we have
\begin{equation}\label{2e3}
\dim I_{L(1,0)} \left(\hspace{-3 pt}\begin{array}{c} L(1,n)\\
L(1, m^{2})\,L(1, n)\end{array}\hspace{-3 pt}\right)=1,
\end{equation}
\begin{equation}\label{2e4}
\dim I_{L(1,0)} \left(\hspace{-3 pt}\begin{array}{c} L(1,k)\\
L(1, m^{2})\,L(1, n)\end{array}\hspace{-3 pt}\right)=0,
\end{equation} for $k\in{\mathbb Z}_{+}$ such that  $k\neq n$.
\end{theorem}

Recall the Heisenberg vertex operator algebra $M(1)$ constructed
from a $d$-dimensional vector space and its subalgebra $M(1)^+$ from
\cite{FLM}. Based on the classification of irreducible modules for
$M(1)^+$ \cite{DN1}, \cite{DN3}, the fusion rules for $M(1)^+$ have
been obtained in \cite{A1} and \cite{ADL}. Here is the result when
$d=1.$
 \begin{theorem}\label{2t2}
Let $M$,$N$ and $T$ be irreducible $M(1)^{+}$-modules. If
$M=M(1,\la)$ such that $\la\neq 0$, then
 $$\dim I_{M(1)^{+}} \left(\hspace{-3 pt}\begin{array}{c} T\\
M\,N\end{array}\hspace{-3 pt}\right)=0 \ \ {\rm or} \ \ 1$$ and
$$\dim I_{M(1)^{+}} \left(\hspace{-3 pt}\begin{array}{c} T\\
M\,N\end{array}\hspace{-3 pt}\right)=1$$ if and only if $(N,T)$ is
one of the following pairs:
$$
(M(1)^{\pm},M(1,\mu))(\la^{2}=\mu^{2}), \ (M(1,\mu),M(1,\nu)), \
(\nu^{2}=(\la\pm \mu)^{2}),$$$$
(M(1)(\theta)^{\pm},M(1)(\theta)^{\pm}), \
(M(1)(\theta)^{\pm},M(1)(\theta)^{\mp}).
$$
\end{theorem}

For the purpose of later discussion we need to study $M(1)^+$ more.
From now on we assume $d=1.$ Recall from \cite{DN1} that
$$J= h(-1)^4{\bf 1} -2h(-3)h(-1){\bf 1} + \frac{3}{2}h(-2)^2{\bf
1}\in M(1)^{+}$$ is a primary vector of weight 4 in $M(1)^{+}.$ Then
$M(1)^+$ is generated by $\omega$ and $J.$ Let $M^{(4)}$ be the
$L(1,0)$-submodule of $M(1)^+$ generated by $J$. Then $J_{7}J=54{\bf
1}.$ Moreover, as a module for $L(1,0),$
$$M(1)^+=\bigoplus_{n\geq 0}L(1,(2n)^2)$$
where $L(1,0)\cong M^{(4)}.$ Following \cite{Z}  we set $$u\ast
v=\Res_z\left( \frac{(1+z)^{\wt(u)}}{z^{1+n}}Y(u,z)v \right)$$ for
homogeneous $u,v\in W$ where $W$ is any vertex operator algebra.
Then
$$J\ast
J=\sum\limits_{j=0}^{4}\left(\begin{array}{c}4\\j\end{array}\right)J_{j-1}J=u^{(0)}+v^{(0)},$$
where $u^{(0)}\in L(1,0)$ and $v^{(0)}\in M^{(4)}$.
\begin{lem}\label{JJ} We have
$$u^{(0)}\in p(\omega)+O(L(1,0)),  \ v^{(0)}\in q(\omega)J+(L(-1)+L(0))M^{(4)}$$
where
$$
p(x)=x(\frac{1816}{35}x^{3}-\frac{212}{5}x^{2}+\frac{89}{10}x-\frac{27}{70}),
$$
$$
q(x)=-\frac{314}{35}x^{2}+\frac{89}{14}x-\frac{27}{70}
$$
and the product is  $\ast.$
\end{lem}

\pf By \cite{DN1}, we have in $M(1)^{+}$
$$
u^{(0)}+v^{(0)}\equiv  p(\omega)+  q(\omega)J\ \mod\,O(M(1)^{+}).$$
Since $J*J\in \sum_{i=0}^8M(1)^+_i$ we see that
$$u^{(0)}\in p_1(\omega)+O(L(1,0)),  \ v^{(0)}\in q_1(\omega)J + (L(-1)+L(0))M^{(4)}$$
where $p_{1}(x)$ is a polynomial with degree $\leq 4$ and $q_{1}(x)$
is a polynomial with degree $\leq 2$. Note that $O(L(1,0)),
(L(-1)+L(0))M^{(4)} \subseteq O(M(1)^{+})$. As in \cite{DN1} we
apply the identity $J*J=u^{(0)}+v^{(0)}$ to the irreducible
$A(M(1)^+)$-modules to conclude that $p_{1}(x)=p(x)$ and
$q(x)=q_{1}(x),$ as desired. \qed

\begin{lem}\label{2l3}  $p(x)$ has no non-zero integer roots.
\end{lem}
\pf It is easy to check that
 $$x_{1}=0, \ x_{2}=\frac{1}{4}, \
 x_{3}=\frac{515+\sqrt{167161}}{1816}, \
 x_{4}=\frac{515-\sqrt{167161}}{1816}$$
 are all the roots of $p(x)$.
 \qed

The following lemma will be used later.
\begin{lem}\label{app1} In $M(1)^{+}$, we have
\begin{align*}
& J_{3}J=-72L(-4){\bf 1}+336L(-2)^2{\bf 1}+\lambda J,\\
& J_{2}J=u^1+\lambda \frac{1}{2}L(-1)J,\\
&
J_{1}J=u^2+\lambda(\frac{28}{75}L(-2)J+\frac{23}{300}L(-1)^{2}J),\\
& J_{0}J=u^3+\lambda(\frac{14}{75}L(-3)J+\frac{14}{75}L(-2)L(-1)J
-\frac{1}{300}L(-1)^{3}J),
\end{align*}
for some $u^i\in L(1,0), i=1,2,3$, $0\neq \lambda\in\C$.
\end{lem}

\pf We first deal with $J_3J.$ Note that $J_iJ\in L(1,0)\oplus
M^{(4)}$ for $i\geq 0.$  Then there exist
$\lambda_1,\lambda_2,\lambda\in \C$ such that
$J_3J=\lambda_1L(-4){\bf 1}+\lambda_2L(-2)^2{\bf 1}+\lambda J.$
Using the commutator formula
$$[L(m), J_n]=[3(m+1)-n]J_{m+m}$$
for $m,n\in\Z$ one can check that
$$(L(-4){\bf 1}, J_3J)=({\bf 1}, L(4)J_3J)=12\times 54,$$
$$(L(-2)^2{\bf 1}, J_3J)=({\bf 1},L(2)^2J_3J)=24\times 54.$$
On the other hand,
$$(L(-4){\bf 1},J_3J)=(L(-4){\bf 1},\lambda_1L(-4){\bf 1}+\lambda_2L(-2)^2{\bf 1})=5\lambda_1+3\lambda_2,$$
$$(L(-2)^2{\bf 1},J_3J)=(L(-2)^2{\bf 1},\lambda_1L(-4){\bf 1}+\lambda_2L(-2)^2{\bf 1})=3\lambda_1+\frac{2}{4}\lambda_2.$$
This implies that $\lambda_1=-72$ and $\lambda_2=336.$ It follows
from the proof of Theorem 4.9 of \cite{DJ2} that $\lambda \ne 0.$
One can also verify $J_{3}J$ directly  with a long computation.

We now prove other relations. We may assume that
\begin{align*}
& J_{2}J=u^1+\mu_{0}L(-1)J,\\
& J_{1}J=u^2+\mu_{1}L(-2)J+\mu_{2}L(-1)^{2}J,\\
& J_{0}J=u^3+\mu_{3}L(-3)J+\mu_{4}L(-2)L(-1)J+\mu_{5}L(-1)^3J
\end{align*}
where $u^i\in L(1,0),i=1,2,3$, $\mu_{j}\in\C,j=0,1,\cdots,5$. Then
$$(J_{2}J, L(-1)J)=(L(1)J_{2}J,J)=4\lambda(J,J),$$
$$(J_{1}J,L(-2)J)=(L(2)J_{1}J,J)=8\lambda(J,J),$$
$$
(J_{1}J,L(-1)^2J)=(L(1)^2J_{1}J,J)=20\lambda(J,J).$$ We also have
$$
(L(-1)J,L(-1)J)=(L(1)L(-1)J,J)=8(J,J),$$
$$
(L(-2)J,L(-2)J)=(L(2)L(-2)J,J)=\frac{33}{2}(J,J),$$
$$
(L(-2)J,L(-1)^2J)=(L(1)^2L(-2)J,J)=24(J,J),$$
$$
(L(-1)^2J,L(-1)^2J)=(L(1)^2L(-1)^2J,J)=144(J,J).$$ These relations
yield the following linear equations
$$
\mu_{0}=\frac{1}{2}\lambda
$$
and
$$\left\{\begin{array}{l}33\mu_{1}+48\mu_{2}=16\lambda\\
24\mu_{1}+144\mu_{2}=20\lambda\end{array}\right.
$$
with solutions
$$\mu_{1}=\frac{28}{75}\lambda, \ \mu_{2}=\frac{23}{300}\lambda.$$

Similarly,
$$
(J_0J,L(-3)J)=(L(3)J_{0}J,J)=12\lambda(J,J),$$
$$
(J_{0}J,L(-2)L(-1)J)=(L(1)L(2)J_{0}J,J)=36\lambda(J,J),$$
$$
(J_{0}J,L(-1)^3J)=(L(1)^3J_{0}J,J)=120\lambda(J,J),$$
$$
(L(-3)J,L(-3)J)=26(J,J), \ (L(-3)J,L(-2)L(-1)J)=40(J,J),$$
$$
(L(-3)J,L(-1)^3J)=96(J,J), \ (L(-2)L(-1)J,L(-2)L(-1)J)=164(J,J), $$
$$
(L(-2)L(-1)J,L(-1)^3J)=624(J,J), \ (L(-1)^3J,L(-1)^3J)=4320(J,J).$$
Then we get the following linear system
$$
\left\{\begin{array}{l} 26\mu_{3}+40\mu_{4}+96\mu_{5}=12\lambda\\
40\mu_{3}+164\mu_{4}+624\mu_{5}=36\lambda\\
96\mu_{3}+624\mu_{4}+4320\mu_{5}=120\lambda\end{array}\right.
$$
with solutions
$$
\mu_{3}=\mu_{4}=\frac{14}{75}\lambda, \
\mu_{5}=-\frac{1}{300}\lambda.$$ The lemma follows.\qed

\section{Search for vector $J$: I}
\def\theequation{3.\arabic{equation}}
\setcounter{equation}{0} In the following discussion throughout  the
paper, we always assume that $V$ is a simple  rational and
$C_{2}$-cofinite vertex operator algebra of central charge 1
satisfying the following conditions:

(1) $V$ is a sum of highest weight modules of $L(1,0)$.

(2) $V=\oplus_{n=0}^{\infty}V_{n}, \ V_{0}=\C {\bf 1}, \ V_{1}=0, \
\dim V_{2}=\dim V_{3}=1, \ \dim V_{4}\geq 3$.

(3) The weights of all the primary vectors in $V$ are perfect
squares.

\begin{remark} As we mentioned in the introduction already, we have dealt with the case  that there exists at
least one primary vector in $V$ whose weight is not a perfect square
in \cite{DJ2}.
\end{remark}

In this section and the next  we look for a primary vector $J'$ of
weight 4 in $V$ such that $J'$ satisfies all relations given in
Lemma \ref{app1} for $J.$ This will help us to show that the vertex
operator subalgebra $U$ generated by $\omega$ and $J'$ is isomorphic
to $M(1)^+$ with identifying $J$ with $J'.$ It turns out that
finding such $J'$ is highly nontrivial and an explicit construction
of intertwining operators for $L(1,0)$ involving modules $L(1,4)$
and $L(1,0)$ plays crucial role in the proof.

Let $X^{1}$ and $X^{2}$ be two subsets of $V$. Set
$$
X^{1}\cdot X^{2}=span\{x_{n}y|x\in X^{1}, y\in X^{2}, n\in{\mathbb
Z}\}.$$ We have the following lemmas from \cite{DJ2} (see also
\cite{DJ1}).
\begin{lem}\label{l1} $V$ is a completely reducible
module for the Virasoro algebra $L(1,0)$.
\end{lem}

\begin{lem}\label{ll2.4}
 Let $u^{1},u^{2}\in V$ be two primary vectors. Let  $U^{1}$ and $U^{2}$ be two $L(1,0)$-submodules of $V$ generated
by $u^{1}$ and $u^{2}$  respectively.  Then
$$
U^{1}\cdot U^{2}=span\{L(-m_{1})\cdots L(-m_{s})u^{1}_{n}u^{2} |
m_{1},\cdots, m_{s}\in {\mathbb Z}_{+}, n\in{\mathbb Z}\}.$$
\end{lem}

For $m\geq 1$, set
$$A_{m^2}=\{v\in V_{m^{2}}| L(n)v=0, n\geq 1\}.$$
Elements in $A_{m^{2}}$ are called primary vectors. Since $V_{1}=0$,
it follows that $m\geq 2$. It is obvious that $V$ carries a
non-degenerate symmetric bilinear form $(\cdot,\cdot)$ such that
$({\bf 1},{\bf 1})=1$ (\cite{FHL}, \cite{L1}). By the assumption
(2), $\dim V_{2}=\dim V_{3}=1$ and  $A_{4}\neq 0$. Let $J'$ be a
non-zero primary vector of weight 4. We may assume that
\begin{equation}\label{eqe4.1}
 (J',J')=54. \ \
\end{equation}
By Theorem \ref{2t1} and Lemma \ref{app1} there exists a primary
vector $u$ of weight $4$  such that
$$J'_{3}J'=-72L(-4){\bf 1}+336L(-2)^2{\bf 1}+27u.$$
It is possible that
$$u=\frac{1}{27}J'_{3}J'+\frac{8}{3}L(-4){\bf
1}-\frac{112}{9}L(-2)^{2}{\bf 1}$$ is zero.

Here is the main result in this section.
\begin{prop}\label{ll4}
The $u$ is a nonzero primary vector of weight 4. In particular,
$J'_{3}J'$ does not lie in $L(1,0).$
\end{prop}

It is not easy to prove this result. We need several lemmas. Let
$V^{(4)}$ be the $L(1,0)$-submodule of $V$ generated by $J'.$ Then
$V^{(4)}$ is isomorphic to $L(1,4).$
\begin{lem}\label{2l1} Let $U$ be the vertex operator subalgebra  of
$V$ generated by $\omega$ and $J'$. If $u\in \C J'$, then $U$ is
linearly spanned by
$$L(-m_{1})\cdots L(-m_{s})J'_{-n_{1}}\cdots J'_{-n_{t}}{\bf 1}
$$
where $m_{1}\geq m_{2}\geq\cdots\geq m_{s}\geq 2$, $n_{1}\geq
n_{2}\geq \cdots\geq n_{t}\geq 1$, $s,t\geq 0$.
\end{lem}
\pf First note that the subspace of $U$ linearly spanned by
$J'_{n_{1}}\cdots J'_{n_{t}}{\bf 1}$ with $n_{i}\in\Z$ is invariant
under the action of $L(m),m\geq -1$. Secondly, we have
\begin{equation}\label{ee10}
[J'_{m},J'_{n}]=\sum\limits_{i=0}^{\infty}\left(\begin{array}{c}m\\i\end{array}\right)
(J'_{i}J')_{m+n-i}.
\end{equation} Since $\wt(J')=4$,
we have $\wt(J'_{i}J')=7-i\leq 7$, for $i\geq 0$.

As $u\in \C J',$ $J'_3J\in L(1,0)+M^{(4)}.$  It follows from Theorem
\ref{2t1} and Lemma \ref{ll2.4} that $J'_{i}J'\in L(1,0)\oplus
V^{(4)}$ for $i\geq 0.$ The lemma is clear. \qed

\begin{lem}\label{2l2} Let $U$ be the vertex operator subalgebra  of
$V$ generated by $\omega$ and $J'$.  If $u\in \C J'$, then the Zhu
algebra $A(U)$ is linearly spanned by
$$\{[\omega]^{s}\ast [J']^t, \ s,t\geq 0\}.$$
\end{lem}
\pf The proof of Theorem 3.5 in \cite{DN1} for the spanning set of
$A(M((1)^+)$ works here.  \qed

\begin{lem}\label{laa1} If $u=0$, then
$$J'\ast J'\equiv p(\omega)+O(L(1,0))\equiv p(\omega)+O(V).$$
\end{lem}

\pf By Theorem \ref{2t1} and Lemma \ref{2l1} we see that $J'\ast
J'\in L(1,0).$ Since $J_{7}J=54{\bf 1}$ and $J'_{7}J'=54{\bf 1}$, it
follows from Theorem \ref{2t1} again that if $u=0$ then
$$J'\ast J'=\sum\limits_{j=0}^{4}\left(\begin{array}{c}4\\j\end{array}\right)J'_{j-1}J'=u^{(0)},$$ where $u^{(0)}$ is the same as in Section 2.
 Then the lemma follows from Lemma \ref{JJ}. \qed

\begin{lem}\label{2l4} Let $U$ be the vertex operator subalgebra  of
$V$ generated by $\omega$ and $J'$.  If $u=0$, then $U=V$.
\end{lem}
\pf By Lemma \ref{2l2}, $A(U)$ is commutative. Suppose that $U\neq
V$. Let $m\in\Z_{+}$ be the smallest positive integer such that
$A_{m^2}\nsubseteq U$. Then $m\geq 2$. Note that $V/U$ is a
$U$-module with the minimal weight $m^2.$ Let $0\neq u^{(m)}\in
A_{m^2}$, $u^{(m)}\notin U$. Then $u^{(m)}+U$ generates a
$U$-submodule of $V/U$ and  let $W$ be the irreducible quotient.  We
denote the image of
 $u^{(m)}+U$  by $v^{(m)}.$  Then $W$
has lowest weight $m^{2}$ and
\begin{equation}\label{ee3} J'_{i}v^{(m)}=0
\end{equation}
for $i\geq 4.$  By Lemmas \ref{2l3} and \ref{laa1} we know that
$J'_3J'_3v^{(m)}=p(L(0))v^{(m)}\ne 0.$ As a result,
\begin{equation}\label{ee6}J'_{3}v^{(m)}\neq 0.
\end{equation}

Since $W$ is completely reducible as an $L(1,0)$-module, $v^{(m)}$
generates an irreducible highest weight $L(1,0)$-submodule
$W^{(m^2)}$ of $W$ isomorphic to $L(1,m^{2})$ with highest weight
vector $v^{(m)}$.  By the skew symmetry, we have $$
J'_{-2}J'=\frac{1}{2}\sum\limits_{i=1}^{\infty}(-1)^{i+1}\frac{L(-1)^{i}}{i!}J'_{-2+i}J'.$$
Then by Theorem \ref{2t1}, Lemma \ref{2l1}, (\ref{ee3}) and the
assumption that $u=0$, we see that
\begin{equation}\label{ee5}V^{(4)}\cdot V^{(4)}\cong L(1,0)\oplus a_{4}L(1, 16).
\end{equation}
 \begin{equation}\label{ee4}V^{(4)}\cdot W^{(m^2)}\cong
 W^{(m^2)}\oplus b_{m+1} L(1,(m+1)^{2})\oplus
b_{m+2} L(1,(m+2)^{2}),\end{equation} where
$a_{4},b_{m+1},b_{m+2}\in\Z$ are nonnegative.

Let ${\mathcal P}$ be the projection from $V^{(4)}\cdot W^{(m^2)}$
to $W^{(m^{2})}$, then ${\mathcal I}(u,z)v={\mathcal P}\cdot
Y(u,z)v$ for $u\in V^{(4)}$, $v\in W^{(m^2)}$ is an intertwining
operator of type
$$\left(\hspace{-3 pt}\begin{array}{c}
L(1,m^2)\\
V^{(4)}\, \ \  L(1,m^2)\end{array}\hspace{-3 pt}\right).$$

 Let $(V_{L}^{+}, {Y}(\cdot,z))$ be the
rank one rational vertex operator algebra with $L=\Z\be$ such that
$(\be,\be)=2m^{2}$.  Set
$$E^{(m)}=e^{\be}+e^{-\be}\in V_{L}^{+}.$$
Then
$$J= h(-1)^4{\bf 1} -2h(-3)h(-1){\bf 1} +
\frac{3}{2}h(-2)^2{\bf 1}\in M(1)^+\subset V_{L}^{+},$$ where
$h=\frac{1}{\sqrt{2}m}\be$. Let $M^{(m^2)}$ be the irreducible
$L(1,0)$-modules in $V_{L}^{+}$ generated by $E^{(m)}.$  From the
construction of $V_{L}^{+}$ we know
\begin{equation}\label{ee7}
M^{(4)}\cdot M^{(m^2)}\cong M^{(m^2)}\oplus L(1,(m+1)^2)\oplus
L(1,(m+2)^{2}),
\end{equation}
\begin{equation}\label{eaa1}
M^{(4)}\cdot M^{(4)}\cong L(1,0)\oplus M^{(4)}\oplus L(1,16).
\end{equation}
Let ${\mathcal Q}$ be the projection from $M^{(4)}\cdot M^{(m^2)}$
to $M^{(m^2)}$. Then ${\mathcal I}'(u,z)v={\mathcal Q}\cdot Y(u,z)v$
for $u\in M^{(4)}$, $v\in M^{(m^2)}$ is an intertwining operator of
type
$$\left(\hspace{-3 pt}\begin{array}{c}
M^{(m^2)}\\
M^{(4)}\, \ \  M^{(m^2)}\end{array}\hspace{-3 pt}\right).$$

 Let $\sigma$ be the $L(1,0)$-module
isomorphism from $M^{(4)}\oplus M^{(m^2)}$ to $V^{(4)}\oplus
W^{(m^2)}$ such that
$$ \sigma(J)=J', \ \sigma(E^{(m)})=v^{(m)}.$$
By Theorem \ref{2t1}, for $u\in M^{(4)}, v\in M^{(m^2)}$,
$$
{\mathcal I}(\sigma u,z)(\sigma v)=c\sigma ({\mathcal I}(u,z)v),$$
for some $c\in {\mathbb C}$. By (\ref{ee6}), $c\neq 0$.  Note that
$$(J'_{7}J')_{-1}v^{(m)}=\sum\limits_{i=0}^{\infty}(-1)^{i}\left(\begin{array}{c}7\\i\end{array}\right)
(J'_{7-i}J'_{-1+i}+J'_{6-i}J'_{i})v^{(m)}.$$ By (\ref{ee4}), for
$i\geq 0$, we have
$$J'_{-1+i}v^{(m)}, \ J'_{7-i}J'_{-1+i}v^{(m)}, \
J'_{i}v^{(m)}, \ J'_{6-i}J'_{i}v^{(m)}\in W^{(m^2)}.$$

On the other hand, we have
$$(J_{7}J)_{-1}E^{(m)}=\sum\limits_{i=0}^{\infty}(-1)^{i}\left(\begin{array}{c}7\\i\end{array}\right)
(J_{7-i}J_{-1+i}+J_{6-i}J_{i})E^{(m)}.$$ By (\ref{ee7}), for $i\geq
0$,
$$J_{-1+i}E^{(m)}, \ J_{7-i}J_{-1+i}E^{(m)}, \
J_{i}E^{(m)}, \ J_{6-i}J_{i}E^{(m)}\in M^{(m^2)}.$$ Thus we have
$$
(J'_{7}J')_{-1}v^{(m)}=c^{2}\sigma((J_{7}J)_{-1}E^{(m)}).$$ Note
that
$$
(J'_{7}J')_{-1}v^{(m)}=54v^{(m)}=\sigma((J_{7}J)_{-1}E^{(m)}).$$ We
deduce that $c^{2}=1.$ Then we may assume that $c=1.$ If $c=-1$,
 replace $J'$ by $-J'.$

 Since
$$
(J'_{3}J')_{3}v^{(m)}=\sum\limits_{i=0}^{\infty}(-1)^{i}\left(\begin{array}{c}3\\i\end{array}\right)
(J'_{3-i}J'_{3+i}+J'_{6-i}J'_{i})v^{(m)},$$
$$
(J_{3}J)_{3}E^{(m)}=\sum\limits_{i=0}^{\infty}(-1)^{i}\left(\begin{array}{c}3\\i\end{array}\right)
(J_{3-i}J_{3+i}+J_{6-i}J_{i})E^{(m)},$$ it follows that
\begin{equation}\label{ee9}
(J'_{3}J')_{3}v^{(m)}=\sigma((J_{3}J)_{3}E^{(m)}).
\end{equation}
Recall from Lemma \ref{app1} that
$$J_{3}J=x+\lambda J$$
for some  $x\in L(1,0)$ and nonzero $\lambda\in \C.$  Suppose that
$$J'_{3}J'=x'+y',$$
where $x'\in L(1,0)$, $y'\in V^{(4)}$. Then by the fact that
$(J,J)=(J',J')$, we have
$$x'=\sigma(x).$$
Then by (\ref{ee9}), for any $w\in M^{(m^2)}$,
$$y'_{3}(\sigma(w))=\sigma(y_{3}w).$$
A straightforward computation shows that $J_3E^{(m)}\ne 0.$ This
implies that $y'\neq 0,$ a contradiction with (\ref{ee5}). This
proves that $U=V.$ \qed

 We are now in a position to prove Proposition \ref{ll4}.

 \pf Suppose that $u=0$, then $U=V$ by Lemma \ref{2l4} and
$$V^{(4)}\cdot V^{(4)}\cong L(1,0)\oplus a_{4}L(1, 16),
$$
where $a_{4}\in\N$. If $a_{4}=0$, then
$$U=L(1,0)\oplus V^{(4)}.$$

Let $W$ be a module for the Virasoro algebra with central charge $c$
such that $W=\bigoplus_{n\in\C}W_n$ where $W_n$ is the eigenspace
for $L(0)$ with eigenvalue $n$ and is finite-dimensional. We define
the $q$-graded dimension of $W$ as
$$\dim_qW=q^{-c/24}\sum_{n\in\C}(\dim W_n)q^n.$$ Denote by $L(c,h)$
the unique irreducible highest weight module for the Virasoro
algebra with central charge $c\in\C$ and highest weight $h\in\C.$
Then
$$\dim_qL(1,h)=\left\{\begin{array}{ll}
{1\over \eta(q)}(q^{n^2/4}-q^{(n+2)^2/4}), & {\rm if} \
h=\frac{1}{4}n^2,n\in\Z\\
\frac{1}{\eta(q)}q^h, & {\rm otherwise.}
\end{array}\right.$$
(cf. \cite{KR}) where
$$\eta(q)=q^{1/24}\prod_{n=1}^{\infty}(1-q^n).$$
Then
$$Z_V(\tau)=ch_qL(1,0)+ch_qL(1,4)=\frac{1-q+q^4-q^9}{\eta(q)}$$
where $q=e^{2\pi i\tau}$ and $\tau$ is a complex variable in the
upper half plane. We sometimes abuse the notation and also denote
$\eta(q)$ by $\eta(\tau).$ Since $V$ is rational and $C_2$-cofinite
we use  the modular invariance result given in \cite{Z} to assert
that  $Z_V(\frac{-1}{\tau})$ should have a $q$-expansion. It is well
known that $\eta(\frac{-1}{\tau})=(-i\tau)^{\frac{1}{2}}\eta(\tau).$
Thus
$$Z_V(\frac{-1}{\tau})=\frac{1-e^{-2\pi i\frac{1}{\tau}}+e^{-2\pi i\frac{4}{\tau}}-e^{-2\pi i\frac{9}{\tau}}}{(-i\tau)^{\frac{1}{2}}\eta(\tau)}$$
which clearly does not have a $q$-expansion. This gives a contradiction.

 So there exists a non-zero
primary vector $u^{(4)}$ of weight $16$ such that
$$u^{(4)}=a_{1} J'_{-9}J'+ x,$$
for some $0\neq a_{1}\in\C$, $x\in L(1,0)\oplus V^{(4)}$. Then by
Lemma \ref{2l1} and (\ref{ee10}), $U$ is linearly spanned by
\begin{equation}\label{ee11}L(-m_{1})\cdots L(-m_{s}){\bf 1}, \ L(-p_{1})\cdots
L(-p_{s}) J'_{-n_{1}}\cdots J'_{-n_{t}}J'_{-9}J', \ L(-p_{1})\cdots
L(-p_{s})J'
\end{equation}
where $m_{1}\geq m_{2}\geq\cdots\geq m_{s}\geq 2$, $n_{1}\geq
n_{2}\geq \cdots\geq n_{t}\geq 9$, $p_{1}\geq p_{2}\geq \cdots \geq
p_{s}\geq 1$. By Theorem \ref{2t1} and (\ref{ee10}), in (\ref{ee11})
we may assume that $n_{t}\geq 17$.

It is easy to see from (\ref{ee11}) that there is no non-zero
primary vector of weight $25.$ If there is no non-zero primary
vector of weight $36$, then by Theorem \ref{2t1},
$J'_{n}J'_{-9}J'\in L(1,0)\oplus V^{(4)}\oplus V^{(16)}$ for
$n\in\Z$ where $V^{(16)}$ is the $L(1,0)$-module generated by
$u^{(4)}.$  This forces that
$$
V\cong L(1,0)\oplus V^{(4)}\oplus V^{(16)}.
$$
The same proof as above gives a contradiction. So there exists a
non-zero primary vector $u^{(6)}$ of weight $36$ such that
$$
u^{(6)}=a_{2} J'_{-17}J'_{-9}J'+x,$$ for some $0\neq a_{2}\in\C$,
$x\in L(1,0)\oplus V^{(4)}\oplus V^{(16)}$. Continuing the process,
we deduce that $V$ is linearly spanned by
 $$L(-m_{1})\cdots L(-m_{s}){\bf 1},\
L(-n_{1})\cdots L(-n_{t})J', \ L(-n_{1})\cdots
L(-n_{t})J'_{-8r-1}\cdots J'_{-9}J', $$ where $m_{1}\geq m_{2}\geq
\cdots m_{s}\geq 2$, $n_{1}\geq n_{2}\geq \cdots n_{t}\geq 1$,
$r\geq 1, s,t\geq 0$ and as a vector space
\begin{equation}\label{ee12}
V\cong \bigoplus _{r=0}^{\infty} L(1,(2r)^{2}). \end{equation} Then
$$Z_V(\tau)=ch_qV=\frac{\sum_{n\geq 0}(-q)^{n^2}}{\eta(q)}=\frac{1}{2\eta(q)}+\frac{\theta_{0,1}(1,q)}{2\eta(q)}$$
where
$$\theta_{0,1}(1,q)=\sum_{n\in\Z}(-1)^nq^{n^2}.$$
is the theta function. It is well known that
$\frac{\theta_{0,1}(1,q)}{2\eta(q)}$ is a modular function over a
congruence subgroup of $SL(2,\Z)$ and $\frac{1}{2\eta(q)}$ is a
modular form of weight $-\frac12.$ On the other hand, $Z_V(\tau)$ is
a component of a vector-valued modular function (cf. \cite{Z},
\cite{KM}, \cite{DM1}). This implies that $\frac{1}{2\eta(q)}$ is a
component of vector-valued modular function over a congruence
subgroup of $SL(2,\Z).$ This is obviously impossible. So $V$ can not
be the form of (\ref{ee12}). \qed

\section{Search for vector $J$: II}
\def\theequation{4.\arabic{equation}}
\setcounter{equation}{0}

In this section we prove that there exists a non-zero primary vector
$X$ of weight $4$ such that $(X,X)\neq 0$ and
$$
X_{3}X=v+c X,$$ for some  $v\in L(1,0)$ and $0\neq  c\in \C.$ Recall
that $A_4$ is the space of primary vectors in $V_4.$ If $\dim
A_{4}=1$, then by Proposition \ref{ll4}, the $J'$ given in Section 3
is the desired element. From now on we assume that $\dim A_{4}\geq
2.$ We will prove that the result is true if $\dim A_4=2$ and then
show that $\dim A_4$ must be $2.$

Assume that $\dim A_{4}=2$. Clearly, there exists $K\in A_{4}$ such
that
\begin{equation}\label{4e1}
(K,K)=54, \ (J',K)=0.
\end{equation}
It follows from Lemma \ref{app1} that  the projection of $J_3J$ to
$L(1,0)$ is
\begin{equation}\label{4ee8}X^0=-72L(-4){\bf 1}+336L(-2)^{2}{\bf 1}.
\end{equation}
Then by the fusion rules of $L(1,0)$ (see Theorem \ref{2t1} and
(\ref{4e1})) we have
$$
J'_{3}J'=X^0+a_{1}J'+b_{1}K, \ K_{3}K=X^0+a_{2}J'+b_{2}K
$$
for some $a_{i},b_{i}\in\C$ with $i=1,2$. If $b_{1}=0$ or $a_{2}=0$,
then by Lemma \ref{ll4}, either $J'$ or $K$ is the desired element
$X\in A_{4}$. So in the following discussion we assume that
$b_{1}\neq 0$, $a_{2}\neq 0$.

From (\ref{4e1}) and Theorem \ref{2t1} we see that $K_iJ'=0$ for all
$i>3.$ Using the skew-symmetry yields $J'_{3}K=K_{3}J'.$ Since
$$
(J'_{3}K,J')=(K,J'_{3}J'), \ (J'_{3}K,K)=(J',K_{3}K)$$ we see that
\begin{equation}\label{eaaa1}
J'_{3}K=K_{3}J'=b_{1}J'+a_{2}K.
\end{equation}
\begin{lem}\label{l5.1} If $\dim A_{4}=2,$ there exists $X\in A_{4}$ such that
\begin{equation}\label{4e2}
X_{3}X=c_{1}X^0+c_{2}X, \end{equation}
 for some
$c_{1},c_{2}\in\C$, where $X^0$ is defined as (\ref{4ee8}).
\end{lem}
\pf For $\mu_{1},\mu_{2}\in \C$,  we have
\begin{align*}
& (\mu_{1}J'+\mu_{2}K)_{3}(\mu_{1}J'+\mu_{2}K)\\
=&
(\mu_{1}^{2}+\mu_{2}^{2})X^0+(\mu_{1}^{2}a_{1}+2\mu_{1}\mu_{2}b_{1}+\mu_{2}^{2}a_{2})J'+(\mu_{1}^{2}b_{1}+2\mu_{1}\mu_{2}a_{2}+\mu_{2}^{2}b_{2})K.
\end{align*}
By the assumption that $b_{1}\neq 0$ and $a_{2}\neq 0$, we may
assume that $\mu_{1}\neq 0$. Then $X=\mu_{1}J'+\mu_{2}K$ satisfies
(\ref{4e2}) for some $c_{1},c_{2}\in\C$ if and only if $\mu_{1}$ and
$\mu_{2}$ satisfy
$$
\dfrac{\mu_{1}^{2}a_{1}+2\mu_{1}\mu_{2}b_{1}+\mu_{2}^{2}a_{2}}{\mu_{1}}=\dfrac{\mu_{1}^{2}b_{1}+2\mu_{1}\mu_{2}a_{2}+\mu_{2}^{2}b_{2}}{\mu_{2}}.
$$
That is,
\begin{equation}\label{4e7}
a_{2}(\frac{\mu_{2}}{\mu_{1}})^{3}+(2b_{1}-b_{2})(\frac{\mu_{2}}{\mu_{1}})^{2}+(a_{1}-2a_{2})\frac{\mu_{2}}{\mu_{1}}-b_{1}=0.
\end{equation}
It is clear that the above equation has solution
$\frac{\mu_2}{\mu_1}\in\C$. The lemma follows. \qed

In the following two lemmas we do not need to assume that $A_4$ is
2-dimensional.

\begin{lem}\label{l5.2}
Let $X\in A_{4}$ be such that $X_{3}X=\mu X^0+\nu X$ for some
$\mu,\nu\in\C$. If $\mu=0$, then $\nu=0$.
\end{lem}
\pf Suppose that $\mu=0$, then $X_{3}X=\nu X.$ By Theorem \ref{2t1}
and Lemma \ref{app1} one deduces that
\begin{equation}\label{4e3}
X_{2}X=\frac{1}{2}\nu L(-1)X,
\end{equation}
\begin{equation}\label{4e4}
X_{1}X=\frac{28}{75}\nu L(-2)X+\frac{23}{300}\nu L(-1)^{2}X,
\end{equation}
\begin{equation}\label{4e5}
X_{0}X=\frac{14}{75}\nu L(-3)X+\frac{14}{75}\nu
L(-2)L(-1)X-\frac{1}{300}\nu L(-1)^{3}X.
\end{equation}
Since $\mu=0$, it follows that $(X,X)=0$. Thus $X_iX=0$ for $i\geq
4.$ In particular, $(X_{4}X)_{2}X=0.$

On the other hand by (\ref{4e3})-(\ref{4e5}), we have
\begin{align*}
& (X_{4}X)_{2}X\\
= &
\sum\limits_{i=0}^{4}(-1)^{i}\left(\begin{array}{c}4\\i\end{array}\right)
(X_{4-i}X_{2+i}-X_{6-i}X_{i})X\\
= & -5X_{4}X_{2}X+4X_{5}X_{1}X-X_{6}X_{0}X\\
= & -5\frac{1}{2}\nu X_{4}L(-1)X+4X_{5}(\frac{28}{75}\nu
L(-2)X+\frac{23}{300}\nu L(-1)^{2}X)\\
& -X_{6}(\frac{14}{75}\nu L(-3)X+\frac{14}{75}\nu
L(-2)L(-1)X-\frac{1}{300}\nu L(-1)^{3}X).
\end{align*}
Using the commutator formula
$$[L(m), X_n]=[3(m+1)-n]X_{m+n}$$
for $m,n\in\Z$ and the fact that $X_iX=0$ for $i\geq 4$ we check
that
$$
X_{4}L(-1)X=4X_{3}X, \ X_{5}L(-2)X=8X_{3}X, \
X_{5}L(-1)^{2}X=20X_{3}X,$$
$$
X_{6}L(-3)X=12X_{3}X, \ X_{6}L(-2)L(-1)X=36X_{3}X, \
X_{6}L(-1)^{3}X=120X_{3}X.$$ Then we deduce that
$$
(X_{4}X)_{2}X=-\frac{12}{25}\nu^{2}X.$$ This proves that $\nu=0$.
\qed

\begin{lem}\label{app3} There are no non-zero elements $X^1,X^2\in
A_{2}$ such that
\begin{equation}\label{fj1}
(X^1,X^1)=0, \ X^1_{3}X^1=0,\ X^2_{3}X^1=\mu X^{(0)}
\end{equation}
for some nonzero $\mu\in \C.$
\end{lem}
\pf Suppose that there are non-zero elements $X^1,X^2\in A_{4}$ such
that (\ref{fj1}) holds.  Note that $Y(X^1,z)X^1\ne 0$ \cite{DL}. Let
$N^i$ be the irreducible $L(1,0)$-modules generated by $X^i$ for
$i=1,2.$
 By Theorem \ref{2t1}, Lemma \ref{ll2.4} and  the skew-symmetry,
$N^1\cdot N^1\cong L(1,16).$ In particular,
$$X^1_{-9}X^1\neq 0, \ X^1_{i}X^1=0, \ i\geq -8.$$

Since $X^2_6X^1\in V_1=0$ we see that
\begin{align*}
& (X^2_{6}X^1)_{0}X^1=
\sum\limits_{i=0}^{6}(-1)^{i+1}\left(\begin{array}{c}6\\i\end{array}\right)
X^1_{6-i}X^2_iX^1\\
=&
\sum\limits_{i=0}^{6}\sum\limits_{j=0}^{7-i}\sum\limits_{s=0}^{8}(-1)^{i+j+s+1}\left(\begin{array}{c}6\\i\end{array}\right)
\left(\begin{array}{c}-8\\j\end{array}\right)\left(\begin{array}{c}8\\s\end{array}\right)X^2_{i+j+s}X^1_{6-i-j-s}X^1\\
= & \sum\limits_{i=0}^{6}\left(\begin{array}{c}6\\i\end{array}\right)\left(\begin{array}{c}-8\\7-i\end{array}\right)K_{15}X_{-9}X\\
=& 8X^2_{15}X^1_{-9}X^1\\
 = & 0.
\end{align*}
This shows that $X^2_{15}X^1_{-9}X^1=0.$ Then we have
\begin{align*}
& (X^2_{7}X^1)_{-1}X^1\\
= &
\sum\limits_{i=0}^{7}(-1)^{i}\left(\begin{array}{c}7\\i\end{array}\right)
X^1_{6-i}X^2_{i}X\\
= &
\sum\limits_{i=0}^{7}\sum\limits_{j=0}^{7-i}\sum\limits_{s=0}^{8}(-1)^{i+j+s}
\left(\begin{array}{c}7\\i\end{array}\right)
\left(\begin{array}{c}-8\\j\end{array}\right)\left(\begin{array}{c}8\\s\end{array}\right)X^2_{i+j+s}X^1_{6-i-j-s}X^1\\
= &
-\sum\limits_{i=0}^{7}\left(\begin{array}{c}7\\i\end{array}\right)\left(\begin{array}{c}-8\\7-i\end{array}\right)X^2_{15}X^1_{-9}X^1\\
 = & 0.
\end{align*}

But by (\ref{fj1}) we know that the projection of $Y(u,z)v$ for
$u\in N^2$ and $v\in N^1$ to $L(1,0)$ is a nonzero intertwining
operator of type $\left(\hspace{-3 pt}\begin{array}{c}
L(1,0)\\
L(4,0)\, \ \ L(4,0)\end{array}\hspace{-3 pt}\right).$  In
particular, $X^2_7X^1=(X^2,X^1){\bf 1}$ is nonzero. This gives a
contradiction and the proof is complete. \qed

\begin{lem}\label{l5.3} Assume that $\dim A_{4}=2$, then there exists $X\in A_{4}$ such that
$(X,X)\neq 0 $ and
$$
X_{3}X=c_{1}X^{(0)}+c_{2}X,$$ for some  $ 0\neq c_{1},0\neq c_{2}\in
\C$.
\end{lem}
\pf By lemma \ref{l5.1}, there exists $X=\mu_{1}J'+\mu_{2}K\in
A_{4}$ such that $\mu_{1},\mu_{2}$ satisfy (\ref{4e7}) and
$$X_{3}X=c_{1}X^{(0)}+c_{2}X, $$
 for some
$c_{1},c_{2}\in\C$. Note that (\ref{4e7}) has three solutions. If
for one solution, $c_{1}\neq 0$, by Lemma \ref{ll4}, $c_{2}\neq 0$.
Then the lemma holds.

Suppose that for all the three solutions of (\ref{4e7}), $c_{1}=0$.
Then $c_{2}=0$ by Lemma \ref{l5.2}. Let
$\nu=\dfrac{\mu_{2}}{\mu_{1}}$, then
$$1+\nu^{2}=1,\ a_{1}+2b_{1}\nu+a_{2}\nu^{2}=0, \
b_{1}+2a_{2}\nu+b_{2}\nu^{2}=0.$$

If $\nu=\sqrt{-1}$, then
$$
a_{1}+2\sqrt{-1}b_{1}-a_{2}=0, \ b_{1}+2\sqrt{-1}a_{2}-b_{2}=0.
$$

 If $\nu=-\sqrt{-1}$,
then
$$
a_{1}-2\sqrt{-1}b_{1}-a_{2}=0, \ b_{1}-2\sqrt{-1}a_{2}-b_{2}=0.
$$
So if (\ref{4e7}) has different solutions, then $b_{1}=a_{2}=0$, a
contradiction with the assumption. This deduces  that all the
solutions of (\ref{4e7}) are $\nu_1=\nu_2=\nu_3=\sqrt{-1}$ or
$\nu_1=\nu_2=\nu_3=-\sqrt{-1}$. Without loss of generality, we
assume that $\nu_{1}=\nu_{2}=\nu_{3}=\sqrt{-1}$. Using the relation
between roots and coefficients of the equation (\ref{4e7}) we see
that $-\sqrt{-1}=\frac{b_1}{a_2}.$ Consequently,
$$b_{1}=-a_{2}\sqrt{-1}, \ a_{1}=-a_{2}, \ b_{2}=a_{2}\sqrt{-1}.$$ We
deduce that
$$
J'_{3}J'=X^{0}+a_{1}(J'+\sqrt{-1}K),
$$
$$
K_{3}K=X^{0}-a_{1}(J'+\sqrt{-1}K).
$$
By Proposition \ref{ll4}, $a_{1}\neq0$. We may assume that
$\mu_{1}=1.$ Then $X=J'+\sqrt{-1}K$. Let $K'=J'-\sqrt{-1}K.$ Then we
have from (\ref{eaaa1}) that
$$(X,X)=0, \ X_{3}X=0,\ (K',K')=0, \ K'_{3}K'=4a_{1}X, \
K'_{3}X=2X^0.$$ This contradicts  Lemma \ref{app3}. \qed

We next establish that $\dim A_{4}\geq 2$ implies $\dim A_{4}=2$.
\begin{lem}\label{l5.4} If $\dim A_{4}\geq 2$, then $\dim A_{4}=2.$
\end{lem}
\pf
 Let $X^{1},\cdots,X^{s}$
be a basis of $A_{4}$ such that
\begin{equation}\label{e4.2.1}
(X^{i},X^{j})=2\delta_{ij}, \ i,j=1,2,\cdots,s. \end{equation}
Recall $X^0$ from (\ref{4ee8}). Then from the discussion on $J'_3J'$
we have
$$
X^{i}_{3}X^{j}=\dfrac{1}{27}\delta_{ij}X^{0}+\sum\limits_{k=1}^{s}a_{ij}^{k}X^{k}
\ $$ for some $a_{ij}^{k}\in\C$, $i,j,k=1,2,\cdots,s.$ The invariant
property
$$(X^{i}_{3}X^{j},X^{k})=(X^{j}_{3}X^{i},X^{k})=(X^{j},X^{i}_{3}X^{k})$$
then gives
\begin{equation}\label{e4.2.2}a_{ij}^{k}=a_{ji}^{k}=a_{ik}^{j}
\end{equation}
for $i,j,k=1,2,\cdots,s$. For $1\leq k\leq s$ we define matrix
$A^{(k)}=(a_{ij}^{k})_{i,j=1}^{s}.$

Using the relation $J_4J=216L(-3){\bf 1}$ and (\ref{e4.2.1}) we see
that $X^{i}_{4}X^{j}=\delta_{ij}8L(-3){\bf 1}$ for
$i,j\in\{1,2,\cdots,s\}.$ This implies for any $k$ that
\begin{equation}\label{e4.2.3} (X^{i}_{4}X^{j})_{2}X^{k}=-64\delta_{ij}X^{k}.\end{equation}
By Lemma \ref{app2}
\begin{align*}
 (X^{i}_{4}X^{j})_{2}X^{k}=u-
 2\sum\limits_{r=1}^{s}\sum\limits_{l=1}^{s}a_{jk}^{r}a_{ir}^{l}X^{l}+
 \dfrac{114}{75}\sum\limits_{r=1}^{s}\sum\limits_{l=1}^{s}a_{ik}^{r}a_{jr}^{l}X^{l}+a\delta_{jk}X^{i}+b\delta_{ik}X^{j}
\end{align*}
for some $a,b\in\C$ where
\begin{align*}
u=&\dfrac{1}{2}a_{jk}^{i}L(-1)X_{4}^{i}X^{i}+a_{ik}^{j}(\dfrac{28}{75}L(-2)+\dfrac{11}{30}L(-1)^{2})X^{j}_{5}X^{j}\\
&
-\dfrac{197}{150}a_{jk}^iL(-1)X_{4}^{j}X^{j}+\dfrac{1}{27}(\dfrac{114}{75}a_{ik}^{j}-2a_{jk}^{i})X^{(0)}.
\end{align*}

Applying $L(1)$ to $X^{i}_{4}X^{i}=8L_{-3}{\bf 1}$ produces
$X^{i}_{5}X^{i}=16L_{-2}{\bf 1}$ for $i=1,2,\cdots,s.$  Since
$(X^{i}_{4}X^{j})_{2}X^{k}=-64\delta_{ij}X^{k}\in A_4$ and $u\in
L(1,0)$
 we see that $u=0$ and
$$-2\sum\limits_{r=1}^{s}\sum\limits_{l=1}^{s}a_{jk}^{r}a_{ir}^{l}X^{l}
 +\dfrac{114}{75}\sum\limits_{r=1}^{s}\sum\limits_{l=1}^{s}a_{ik}^{r}a_{jr}^{l}X^{l}+a\delta_{jk}X^{i}+b\delta_{ik}X^{j}=-64\delta_{ij}X^{k}.
$$
Comparing the coefficients of $X^l$ of both sides and varying $i,j$
we have  for all $l,k$  that
\begin{equation}\label{e4.2.4}
-\dfrac{36}{75}A^{(k)}A^{(k)}=-(a+b)E_{kk}-64I.
\end{equation}
\begin{equation}\label{e4.2.5}
\dfrac{114}{75}A^{(k)}A^{(l)}-2A^{(l)}A^{(k)}=-aE_{lk}-bE_{kl}-64\delta_{kl}I.
\end{equation}
\begin{equation}\label{e4.2.6}
\dfrac{114}{75}A^{(l)}A^{(k)}-2A^{(k)}A^{(l)}=-aE_{kl}-bE_{lk}-64\delta_{lk}I,
\end{equation}
where $I$ is the identity matrix and $E_{pq}=(e_{ij})_{i,j=1}^{s}$
such that $e_{ij}=\delta_{ip}\delta_{jq}$. Then we deduce that for
$1\leq k,l\leq s$, $k\neq l$,
\begin{equation}\label{e4.2.7}
\dfrac{11}{75}A^{(k)}A^{(l)}=\dfrac{1}{144}[(25b+19a)E_{lk}+(19b+25a)E_{kl}].
\end{equation}

Now suppose that $s\geq 3$.  For $1\leq k\leq s$, denote by
$r(A^{(k)})$ the rank of $A^{(k)}$. By (\ref{e4.2.7}),
$r(A^{(k)})\leq s-1$. It follows from (\ref{e4.2.4}) that $a+b=-64$
and $ \dfrac{9}{75}A^{(k)}A^{(k)}=16(I-E_{kk}).$ Using
$\dfrac{11\times 9}{75}(A^{(1)}A^{(1)})A^{(2)}=\dfrac{11\times
9}{75}A^{(1)}(A^{(1)}A^{(2)})$ gives
\begin{equation}\label{e4.2.8}
176\left[\begin{array}{cccccc}0 &\   0 & \ \cdots & \  0 &  \ 0\\
a_{21}^{2} & \  a_{22}^{2} & \ \cdots & \  a_{2,s-1}^{2}& \ a_{2s}^{2}\\
\vdots & \ \vdots  &\ \cdots & \  \vdots&\ \vdots\\
a_{s1}^{2} & \ a_{s2}^{2}&\ \cdots & \ a_{s,s-1}^{2}& \ a_{ss}^{2}
\end{array}\right]=\dfrac{1}{16}\left[\begin{array}{cccccc} d_{1}a_{12}^{1}&\   d_{2}a_{11}^{1} & \ 0 & \cdots  &  \ 0\\
d_{1}a_{22}^{1} & \  d_{2}a_{21}^{1} & \ 0 & \cdots & \  0\\
\vdots & \ \vdots  &\ \vdots & \  \cdots&\ \vdots\\
d_{1}a_{s2}^{1} & \ d_{2}a_{s1}^{1}&\ 0 & \cdots & \ 0
\end{array}\right]
\end{equation}
where $d_{1}=25b+19a$ and $ d_{2}=19b+25a$. Similarly, by the fact
that $$ (A^{(1)}A^{(2)})A^{(2)}=A^{(1)}(A^{(2)}A^{(2)}),$$ we have
\begin{equation}\label{e4.2.10}
176\left[\begin{array}{cccccc}
a_{11}^{1} & \ 0 & \  a_{13}^{1} & \ \cdots & \ a_{1s}^{1}\\
a_{21}^{1} & \ 0 & \  a_{23}^{1} & \ \cdots & \ a_{2s}^{1}\\
\vdots & \ \vdots  &\ \vdots & \  \cdots&\ \vdots\\
a_{s1}^{1} & \ 0 & \  a_{s3}^{1}&\ \cdots &  \ a_{ss}^{1}
\end{array}\right]=\dfrac{1}{16}\left[\begin{array}{cccccc} d_{2}a_{21}^{2}&\   d_{2}a_{22}^{2}  & \cdots  &  \ d_{2}a_{s2}^{2}\\
d_{1}a_{11}^{2} & \  d_{1}a_{12}^{2} & \cdots & \  d_{1}a_{1s}^{2}\\
0& \ 0 & \ \cdots & 0\\
\vdots & \ \vdots  & \  \cdots&\ \vdots\\
0& \ 0 & \ \cdots & 0
\end{array}\right].
\end{equation}

By (\ref{e4.2.8}), $a_{ij}^2=0$ for $i\geq 2, j>2.$   Then by
(\ref{e4.2.2}), $a_{ji}^{2}=0$ with $i\geq 2, j>2.$ Using
(\ref{e4.2.8}) again asserts
\begin{equation}\label{e4.2.9} d_{2}a_{j1}^{1}=0, \ j=3,4,\cdots,s.\end{equation}
Assume  $d_{2}\neq 0$ and $d_{1}\neq 0.$ It follows from
(\ref{e4.2.9}) that $a_{j1}^{1}=0$ for $j\geq 3.$ Using
(\ref{e4.2.8}) and  (\ref{e4.2.10})  also gives
$a_{11}^{1}=a_{12}^{1}=0.$ This implies that $a_{11}^k=0$ for all
$k$ and
\begin{equation}\label{e4.2.11}X^{1}_{3}X^{1}=\dfrac{1}{27}X^{(0)}.
\end{equation}
This contradicts  Proposition \ref{ll4}. So $d_{1}=0$ or $d_{2}=0$.
If $d_{1}=25b+19a=0$ then $d_{2}=19b+25a\neq0$ as $a+b=-64.$ We have
$a_{1j}^{1}=0$ for $j=1$ or $j>2.$ By (\ref{e4.2.10}),
$a_{12}^1=a_{21}^1=0$. So we have (\ref{e4.2.11}) again. Similarly
if $d_{2}=0$ then $d_{1}\neq 0,$ $a_{1j}^1=0$ for all $j$ and
(\ref{e4.2.11}) holds. The proof is complete. \qed

\section{The subalgebra $M(1)^+$  of $V$}
\def\theequation{5.\arabic{equation}}
\setcounter{equation}{0}

In this section, we prove that there is a vertex operator subalgebra
$U$ of $V$ isomorphic to $M(1)^{+}$. By Lemmas \ref{l5.3} and
\ref{l5.4}, there exists $J'\in A_{4}$ such that
\begin{equation}\label{e5.1}
J'_{3}J'=X^{0}+c J'
\end{equation}
for some $0\neq c\in \C$. Let $U$ be the subalgebra of $V$ generated
by $\omega$ and $J'$. Recall that $V^{(4)}$ is the irreducible
$L(1,0)$-submodule of $V$ generated by $J'$. From the skew-symmetry
$$
Y(J',z)J'=e^{L(-1)z}Y(J',-z)J'$$ we see that
$$
J'_{-2}J'=-J'_{-2}J'+\sum_{j=1}^{9}(-1)^{j+1}\frac{1}{j!}L(-1)^{j}J'_{-2+j}J'.
$$
This together with Theorem \ref{2t1}  deduces that
\begin{equation}\label{equ1}
V^{(4)}\cdot V^{(4)}\subseteq L(1,0)\oplus V^{(4)}\oplus L(1,16).
\end{equation}
The proof of the following lemma is similar to that of Proposition
\ref{ll4}.
\begin{lem}\label{l6.1}
The vertex operator subalgebra $U$ is not equal to the whole algebra
$V.$
\end{lem}

Recall Lemma \ref{JJ} and $J_3J=X^{(0)}+\lambda J$ from Lemma
\ref{app1}. The following is an analogue to Lemma \ref{JJ}.
\begin{lem}\label{J'} We have
$J'\ast J'=u^{(0)}+\frac{c}{\lambda}v^{(0)}$ in $U$ where
$$u^{(0)}\in p(\omega)+O(L(1,0)),  \ v^{(0)}\in q(\omega)J'+(L(-1)+L(0))V^{(4)}$$
where $p(x)$ and $q(x)$ are defined in Lemma \ref{app1}.
\end{lem}
\pf First we define four projections:
 $$p_1: M(1)^+\to L(1,0),\  p_2: M(1)^+\to M^{(4)},\ p_1': U\to L(1,0),\ p_2': U\to  V^{(4)}.$$
 Then $p_1Y(u,z)v$ and $p_2Y(u,z)v$ for $u,v\in M^{(4)}$ define two intertwining operators of types
$\left(\hspace{-3 pt}\begin{array}{c} L(1,0)\\
M^{(4)}\,M^{(4)}\end{array}\hspace{-3 pt}\right)$ and $\left(\hspace{-3 pt}\begin{array}{c} M^{(4)}\\
M^{(4)}\,M^{(4)}\end{array}\hspace{-3 pt}\right),$ respectively.
Similarly, $p_1'Y(u,z)v$ and $p_2'Y(u,z)v$ for $u,v\in V^{(4)}$
define two intertwining operators of types
$\left(\hspace{-3 pt}\begin{array}{c} L(1,0)\\
V^{(4)}\,V^{(4)}\end{array}\hspace{-3 pt}\right)$ and $\left(\hspace{-3 pt}\begin{array}{c} V^{(4)}\\
V^{(4)}\,V^{(4)}\end{array}\hspace{-3 pt}\right),$ respectively.
Identify $L(1,0)$ in both $M(1)^+$ and $U.$ Since
$J_3J=X^{(0)}+\lambda J$ and $J_3'J'=X^{(0)}+cJ'$, the result
follows from Theorem \ref{2t1} and Lemma \ref{JJ} immediately. \qed

By Lemmas \ref{l6.1} and \ref{l5.4}  there exists $2\leq k\in\Z_{+}$
such that $U_m=V_m$ for $m< k^2$ and $U_{k^2}\ne V_{k^2}.$ Take
$F\in A_{k^2}\notin U.$ As in the proof of Lemma \ref{2l4},  $V/U$
is a $U$-module with the minimal weight $k^2$ and
 $F^{(k)}+U$ generates a
$U$-submodule $W$ of $V/U.$ Let $\bar{W}$ be the irreducible
quotient of $W$.
 We denote the image of $F+U$  in $\bar{W}$ by $a.$  Note that
 $\bar{W}$ has the lowest weight $k^2$ and $a$ generates an
irreducible $L(1,0)$-submodule $W^{(k^2)}$ of $\bar{W}$. By Lemma
\ref{2l2}, $A(U)$ is commutative. It follows that the lowest weight
subspace is one-dimensional. Then
$$J'_{3}a\in\C v, \  J'_{n}a=0,$$
for $n\geq 4$. Since $p(L(0))a\ne 0$ (see Lemma \ref{2l3}) we
immediately deduce from Lemma \ref{J'} that $J'_{3}a\neq 0$.

 Let $V_{L}^{+}$ be the rational vertex operator algebra associated
to the definite even lattice $L=\Z\al$ such that $(\al,\al)=2k^{2}$.
Set
$$E=e^{\al}+e^{-\al}\in V_{L}^{+},$$
$$J= h(-1)^4{\bf 1} -2h(-3)h(-1){\bf 1} +
\frac{3}{2}h(-2)^2{\bf 1}\in V_{L}^{+},$$ where
$h=\frac{1}{\sqrt{2}k}\al$.  We identify the Virasoro vertex
operator subalgebra $L(1,0)$ in $V$ and $V_{L}^{+}$. Let $M^{(k^2)}$
be the irreducible $L(1,0)$-module in $V_{L}^{+}$ generated by $E.$
Then there exists an $L(1,0)$-module isomorphisms $\sigma$ from
$M^{(4)}$ to $V^{(4)}$ and $M^{(k^2)}$ to $W^{(k^2)}$ such that
$\sigma J=J'$ and $\sigma E=a.$

Let $P$ and $P'$ be the projections of $V_{L}^{+}$ to $M^{(k^2)}$
and $\bar{W}$ to $W^{(k^2)}$ respectively. Then $I(u,z)w=P\cdot
Y(u,z)w$ and  $ I'(u',z)w'=P'\cdot Y(u',z)w'$ for $u\in M^{(4)},
w\in M^{(k^2)},$ $ u'\in V^{(4)}, w'\in W^{(k^2)}$ are
intertwining operators of types
$\left(\hspace{-3 pt}\begin{array}{c} M^{(k^2)}\\
M^{(4)}\,M^{(k^2)}\end{array}\hspace{-3 pt}\right)$ and  $\left(\hspace{-3 pt}\begin{array}{c} W^{(k^2)}\\
V^{(4)}\,W^{(k^2)}\end{array}\hspace{-3 pt}\right),$ respectively.
Let $Q$ and $Q'$ be the projections of $V_{L}^{+}$ to $M^{(4)}$ and
$V$ to $V^{(4)},$ respectively. Then $ {\cal I}(u,z)w=Q\cdot
Y(u',z)w'$ and  ${\cal I}'(u',z)w'$ for $u,w\in M^{(4)},$ $u',w'\in
V^{(4)}$ are   intertwining operators of type
$\left(\hspace{-3 pt}\begin{array}{c} M^{(4)}\\
M^{(4)}\,M^{(4)}\end{array}\hspace{-3 pt}\right)$ and $\left(\hspace{-3 pt}\begin{array}{c} V^{(4)}\\
V^{(4)}\,V^{(4)}\end{array}\hspace{-3 pt}\right),$ respectively.

It follows from Theorem \ref{2t1} that $J'_{n}a\in W^{(k^2)}$ and
$$
(J'_{3}J')_{n}a=\sum\limits_{i=0}^{\infty}(-1)^{i}\left(\begin{array}{c}3\\i\end{array}\right)
(J'_{3-i}J'_{n+i}+J'_{3+n-i}J'_{i})a\in W^{(k^2)}
$$
for $n\geq 3-2k$. Using the proof of Lemma 4.5 and Lemma 4.7 in
\cite{DJ2} we have
\begin{lem}\label{l6.2}
$$
I'(\sigma(u),z)\sigma(w)=\sigma(I(u,z)w), \ \ {\cal
I'}(\sigma(u),z)\sigma(v)=\sigma({\cal I}'(u,z)v)
$$
for $u,v\in M^{(4)}, w\in W^{(k^2)}$. In particular, let $\la\in\C$
be such  that
$$
J_{3}J=X^{(0)}+\lambda J,$$ then
$$
J'_{3}J'=X^{(0)}+\lambda J'.$$
\end{lem}

Recall that $J'_{n}w\in W^{(k^2)}$ for $n\geq 3-2k.$ That is,
$J'_nw$ is a linear combination of vectors of form $L(-n_1)\cdots
L(-n_k)w$ with $3-n=-n_1-n_2\-\cdots -n_k.$ Using the relation
$[L(m),J'_n]=(3(m+1)-n)J'_{m+n}$ we see that if $n\geq 15-2k$ then
$$
(J'_{-9}J')_{n}a=\sum\limits_{i=0}^{\infty}(-1)^{i}\left(\begin{array}{c}-9\\i\end{array}\right)
(J'_{-9-i}J'_{n+i}+J'_{-9+n-i}J'_{i})a\in W^{(k^2)}.$$ A proof
similar to that of Lemma 4.10 in \cite{DJ2} gives that
\begin{lem}\label{l6.3}
There  exist a non-zero primary element $u^4$ of weight 16 in
$V_{L}^{+}$ and a non-zero primary element $v^4$ of weight 16 in $U$
such that the isomorphism $\sigma$ from $L(1,0)\bigoplus M^{(4)}$ to
$L(1,0)\bigoplus V^{(4)}$ can be extended to an isomorphism $\sigma$
from $L(1,0)\bigoplus M^{(4)}\bigoplus M^{(16)}$ to $L(1,0)\bigoplus
V^{(4)}\bigoplus V^{(16)}$ such that $\sigma(u^4)=v^4$ and for
$n\in\Z$, $u,v\in L(1,0)\bigoplus M^{(4)}$,
\begin{equation}
(\sigma u)_{n}(\sigma v)=\sigma(u_{n}v),
\end{equation}
where $M^{(4^2)}$ is the irreducible  $L(1,0)$-submodule of
$V_{L}^{+}$ generated by  $u^4$ and $V^{(4^2)}$ is the irreducible
$L(1,0)$-submodule of $U$ generated by $v^4$.  In particular, for
$m,n\in\Z$,
\begin{equation}\label{e6.3}
[J'_{m},J'_{n}]=\sum\limits_{j=0}^{\infty}\left(\begin{array}{c}m\\j\end{array}\right)(\sigma(J_{j}J))_{m+n-j}.
\end{equation}
\end{lem}

We  can now have a ``nicer'' spanning set for $U.$
\begin{lem}\label{l6.4}
$U$ is linearly spanned by
$$L(-m_{s})\cdots L(-m_{1}){\bf 1},\
L(-n_{s})\cdots L(-n_{1})J',\ L(-n_{s})\cdots
L(-n_{1})J'_{-8t-1}\cdots J'_{-17}J'_{-9}J'$$
 where
$m_{s}\geq m_{s-1}\geq \cdots m_{1}\geq 2$, $n_{s}\geq n_{s-1}\geq
\cdots n_{1}\geq 1$, $t\geq 1$, $s\geq 0$.
\end{lem}
\pf By Lemma \ref{2l1} and (\ref{e5.1}) $U$ is linearly spanned by
$$L(-m_{s})\cdots L(-m_{1}){\bf 1}, \ L(-n_{s})\cdots L(-n_{1})J',\ L(-n_{s})\cdots L(-n_{1})J'_{-p_{t}}\cdots
J'_{-p_{1}}J'
$$
where $m_{s}\geq m_{s-1}\geq\cdots\geq m_{1}\geq 2$, $n_{s}\geq
n_{s-1}\geq \cdots\geq n_{1}\geq 1$, $p_{t}\geq p_{t-1}\geq
\cdots\geq p_{1}\geq 1$, $s,t\geq 0$.

It follows from  (\ref{e5.1}) and Theorem \ref{2t1} that
$$L(-n_{s})\cdots L(-n_{1})J'_{-p_{1}}J'\in
L(1,0)\bigoplus V^{(4)}$$ for $p_{1}\leq 8.$ So we can assume that
$p_{1}\geq 9$. Using (\ref{equ1}) gives
$$L(-n_{s})\cdots L(-n_{1})J'_{-p_{1}}J'\in
L(1,0)\bigoplus V^{(4)}\bigoplus V^{(16)}$$ for $p_{1}\geq 9$. By
Lemma \ref{l6.3}, $V^{(4)}\cdot V^{(4)}\cong L(1,0)\bigoplus
V^{(4)}\bigoplus V^{(16)}$ and  there is a non-zero primary vector
$v^{(16)}$ of weight $16$ in $V^{(4)}\cdot V^{(4)}$ such that
$v^{(16)}=x+aJ'_{-9}J'$ for some $x\in L(1,0)\bigoplus V^{(4)}$ and
$0\neq a\in\C$. So we may assume that $p_{1}=9$.

If there exists a non-zero primary vector $u$ of weight $25$ then
$u=u^1+aJ'_{-6}J'_{-9}J'$ for some $u^1\in L(1,0)\bigoplus
V^{(4)}\bigoplus V^{(16)}$ and $0\neq a\in\C$. Note that $J'_{-6}J',
J'_{j}J'\in L(1,0)\bigoplus V^{(4)}$ for $j\geq 0$. So
$$J'_{-6}J'_{-9}J'=J'_{-9}J'_{-6}J'+
\sum\limits_{j=0}^{\infty}\left(\begin{array}{c}-6\\j\end{array}\right)(J'_{j}J')_{-15-j}J'\in
L(1,0)\bigoplus V^{(4)}\bigoplus V^{(16)}.$$ This proves that there
is no non-zero primary vector of weight $25.$  By Theorem \ref{2t1}
again
$$L(-n_{s})\cdots L(-n_{1})J'_{-p_{2}}
J'_{-p_{1}}J'\in L(1,0)\bigoplus V^{(4)}\bigoplus V^{(16)}$$ for
$p_{2}<17$ and
$$L(-n_{s})\cdots L(-n_{1})J'_{-p_{2}}
J'_{-p_{1}}J'\in L(1,0)\bigoplus V^{(4)}\bigoplus V^{(16)}\bigoplus
L(1,36)$$ for $p_{2}\geq 17$. So we may assume that $p_{2}=17.$
Continuing in this way gives the result. \qed

The following is the main result of this section.
\begin{theorem}\label{t6.1}
There is a vertex operator algebra isomorphism $\sigma$ from
$M(1)^{+}$ to $U$ such that $\sigma \omega=\omega$ and
$\sigma(J)=J'.$
\end{theorem}
\pf Let $u^4$, $v^4$ and $\sigma$ be the same as  in Lemma
\ref{l6.3}. Then there exist $x^{(4)}\in L(1,0)\bigoplus M^{(4)} $
and $0\neq a_{1}\in\C$ such that
$$u^4=a_{1}J_{-9}J+x^{(4)}, \
v^4=a_{1}J'_{-9}J'+\sigma(x^{(4)}).
$$
Moreover, $(u,v)=(\sigma(u),\sigma(v))$ for $u,v\in
L(1,0)+M^{(4)}+M^{(16)}.$ From the construction of $M(1)^{+}$
\cite{DG},  there exists a non-zero primary element $u^6$ of weight
36 in $M(1)^{+}$  such that $J_{-17}J_{-9}J=u^6+x^{(6)}, $ where
$x^{(6)}\in L(1,0)\bigoplus M^{(4)}\bigoplus M^{(16)}.$ It is
obvious that
$(J_{-17}J_{-9}J,J_{-17}J_{-9}J)=(u^6,u^6)+(x^{(6)},x^{(6)}).$ Set
$M^{(0)}=V^{(0)}=L(1,0).$ Let $P_i$ and $Q_i$ be the projections of
$V_{L}^{+}$ and $U$ to  $M^{(i)}$ and $V^{(i)},$ respectively for
$i=0,4,16$. Then $I^i(u,z)v=P_iY(u,z)v$ and $J^i(\sigma u,z)\sigma
v=Q_iY(\sigma u,z)\sigma v$ for $u\in M^{(4)}$, $v\in M^{(16)}$ are
intertwining
operators of type $\left(\hspace{-3 pt}\begin{array}{c} M^{(i)}\\
M^{(4)}\,M^{(16)}\end{array}\hspace{-3 pt}\right)$ and $\left(\hspace{-3 pt}\begin{array}{c} V^{(i)}\\
V^{(4)}\,V^{(16)}\end{array}\hspace{-3 pt}\right)$ respectively for
$i=0,4,16.$

By Lemma \ref{l6.3} for $n\geq -16$,
\begin{align*}
&J'_{n}J'_{-9}J'=\sum\limits_{j=0}^{\infty}
\left(\begin{array}{c}n\\j\end{array}\right)(J'_{j}J')_{-9+n-j}J'+J'_{-9}J'_nJ'\\
& =\sum\limits_{j=0}^{\infty}\left(\begin{array}{c}n\\j\end{array}\right)\sigma((J_{j}J)_{-9+n-j}J)+\sigma(J_{-9}J_nJ)\\
&=\sigma(J_{n}J_{-9}J).
\end{align*}
Note that $J_{-9}J\in L(1,0)\bigoplus M^{(4)}\oplus M^{(16)}$ from
the structure of $M(1)^+$ and there exist $m,n\geq 3$ such that the
projections of $J_mJ_{-9}J$ and $J_nJ_{-9}J$ to $M^{(4)}$ and
$M^{(16)}$ are nonzero. We also know that $J'_{-9}J'\in L(1,0)\oplus
V^{(4)}\oplus V^{(16)}$. It follows that
\begin{equation}\label{e5.5}
J^i(\sigma u,z)\sigma v=\sigma (I^i(u,z)v)
\end{equation}
for $i=2, 4$ and $u\in M^{(4)}$ and $v\in M^{(16)}.$ By Theorem
\ref{2t1} and (\ref{e5.5}), there exists a  non-zero primary vector
$v^6$ of weight 36 in  $U$ such that
$J'_{-17}J'_{-9}J'=v^6+\sigma(x^{(6)}).$ By (\ref{e5.5}) again
\begin{align*}
&(J'_{-17}J'_{-9}J',J'_{-17}J'_{-9}J')=(J'_{-9}J',J'_{23}J'_{-17}J'_{-9}J')\\
&=\sum\limits_{j=0}^{\infty}\left(\begin{array}{c}23\\j\end{array}\right)(J'_{-9}J',(J'_{j}J')_{6-j}J'_{-9}J')\\
&=\sum\limits_{j=0}^{\infty}\left(\begin{array}{c}23\\j\end{array}\right)(J_{-9}J,(J_{j}J)_{6-j}J_{-9}J)\\
&=(J_{-17}J_{-9}J,J_{-17}J_{-9}J).
\end{align*}
In particular, $(v^6,v^6)=(u^6,u^6).$

Let $M^{(36)}$  be the irreducible $L(1,0)$-submodule of $V_{L}^{+}$
 generated by $u^6$ and $V^{(36)}$ the irreducible $L(1,0)$-submodule of $U$ generated by $v^6$. Then the
$L(1,0)$-module  isomorphism $\sigma$ from
$M^{(0)}+M^{(4)}+M^{(16)}$ to $V^{(0)}+V^{(4)}+V^{(16)}$ can be
extended to the $L(1,0)$-module isomorphism $\sigma$ from
$M^{(0)}+M^{(4)}+M^{(16)}+M^{(36)}$  to
$V^{(0)}+V^{(4)}+V^{(16)}+V^{(36)}$ such that
$$\sigma(u^6)=v^6,\  \sigma(J_{-17}J_{-9}J)=J'_{-17}J'_{-9}J'$$
and for $n\geq -24$,
$$
\sigma(J_{n}J_{-17}J_{-9}J)=J'_{n}J'_{-17}J'_{-9}J'.$$ Similarly, we
have for any $n\in\Z$,
$$
\sigma(J_{n}J_{-17}J_{-9}J)=J'_{n}J'_{-17}J'_{-9}J'.$$ Continuing
the above process and using Lemma \ref{l6.4}, we deduce that as a
vector space,
$$U\cong \bigoplus_{m=0}^{\infty}L(1, 4m^{2})$$ and there is an $L(1,0)$-module isomorphism $\sigma$ from $M(1)^{+}$ to
$U$ such that
$$\sigma(x)=x, \sigma(J)=J', \ \sigma(J_{-8m-1}\cdots J_{-17}J_{-9}J)=J'_{-8m-1}\cdots J'_{-17}J'_{-9}J'$$
for $x\in L(1,0)$and  $m\geq 1$,  and for $n\in\Z$ $$
\sigma(J_{n}J_{-8m-1}\cdots J_{-17}J_{-9}J)=J'_{n}J'_{-8m-1}\cdots
J'_{-17}J'_{-9}J'.$$ Then it follows from Theorem 5.7.1 in \cite{LL}
that $\sigma$ is an isomorphism of vertex operator algebras. \qed
\section{Identification of  $V$ with $V_L^+$}
\def\theequation{6.\arabic{equation}}
\setcounter{equation}{0}

In this section, we will prove that $V$ is isomorphic to the
rational vertex operator algebra $V_{\Z\al}^{+}$ where
$(\al,\al)=2k^2$ and $k\geq 2$ is the smallest positive integer such
that $U_{k^2}\ne V_{k^2}.$

By Theorem \ref{t6.1}, $U$ is a simple vertex operator subalgebra of
$V$ isomorphic to $M(1)^+.$ Then the restriction of the
non-degenerate invariant symmetric bilinear form on $V$ to $U$ is
non-degenerate. We identify $U$ with $M(1)^+.$ Then
\begin{equation}\label{e6.1}
V=M(1)^{+}\bigoplus (M(1)^{+})^\bot,
\end{equation}
where $(M(1)^{+})^\bot$ is the orthogonal complement of $M(1)^+$ in
$V$ with respect to the bilinear form. Clearly, $(M(1)^{+})^\bot$ is
also an $M(1)^+$-module. Since $V$ is rational, it follows that
$(M(1)^{+})^\bot\neq 0$. Then $(M(1)^{+})^\bot\cong \bigoplus_{m\geq
k}c_{m}L(1,m^2)$ for some $k\geq 2$ such that $c_{k}\neq 0$ and
$c_{m}\in\N$. Then $(M(1)^{+})^\bot=\sum_{m\geq
k}(M(1)^{+})^\bot_m.$ It is obvious that
\begin{equation}\label{e6.2}
L(m)(M(1)^{+})^{\bot}_{k^2}=J'_n(M(1)^{+})^{\bot}_{k^2}=0,
J'_{3}(M(1)^{+})^\bot_{k^2}\subseteq (M(1)^{+})^\bot_{k^2}
\end{equation}
for $m\geq 1$ and $n\geq 4.$

\begin{lem}\label{lem6.1} Let $W$ be an $M(1)^+$-module such that $W$ is a completely reducible $L(1,0)$-module.
Let $v\in W$ be a non-zero primary element of weight $n^2$ ($n\geq
2$ is an integer) such that $ J'_{3}v\in\C v$ and $J'_{m}v=0$ for
$m\geq 4$. Then the $M(1)^{+}$-submodule $N$ of $W$ generated by $v$
is irreducible.
\end{lem}
\pf Since $N$ is a completely reducible $L(1,0)$-module, it follows
from Theorem \ref{2t1} that
$$
N=\bigoplus_{p=0}^{\infty}c_{p}L(1,(n+p)^{2}),$$ for some $c_{p}\in
\N$, $p=0,1,2,\cdots$. Let $\bar{N}$ be the irreducible quotient of
$N$. Then (see \cite{DG} and \cite{DN1})
$$\bar{N}=\bigoplus_{p=0}^{\infty}L(1,(n+p)^{2})$$
as an $L(1,0)$-module.  So $c_{p}\geq 1,p=0,1,2,\cdots$.  Obviously,
$N$ is linearly spanned by
$$L(-m_{s})\cdots L(-m_{1})v, \
L(-n_{s})\cdots L(-n_{1})J'_{-p_{t}}\cdots J'_{-p_{1}}v$$ where
$m_{s}\geq m_{s-1}\geq\cdots\geq m_{1}\geq 1$, $n_{s}\geq
n_{s-1}\geq \cdots\geq n_{1}\geq 1$, $p_{t}\geq p_{t-1}\geq
\cdots\geq p_{1}$, $s,t\geq 0$. By Theorem \ref{2t1} and the fact
that $J'_{3}v\in\C v$, we may assume that $p_{1}\geq 2n-2$ in
(\ref{eq8.3}). Using a similar proof given in Lemma \ref{l6.4} shows
that $N$ is, in fact, spanned by
\begin{equation}\label{eq8.3}
L(-m_{s})\cdots L(-m_{1})v, \  \  L(-n_{s})\cdots
L(-n_{1})J'_{2-2(n+i)}\cdots J'_{2-2(n+1)}J'_{2-2n}v
\end{equation}
 where
$m_{s}\geq m_{s-1}\geq \cdots m_{1}\geq 1$, $n_{s}\geq n_{s-1}\geq
\cdots n_{1}\geq 1$, $i\geq 0$, $s\geq 0$. Let $N^i$ denote the
subspace of $N$ spanned by the elements given in (\ref{eq8.3}) for
fixed $i.$ Then each $N^i$ is an $L(1,0)$-submodule of $N$ and
$\sum_{i}N_i=N.$

Let $v^1$ be a highest weight vector of weight $(n+1)^2$ in $N$.
Then by (\ref{eq8.3})
$$v^1=u^1+aJ'_{2-2n}v$$ for some $u^1\in L(1,n^2), 0\neq
a\in\C$. Thus $N^0\cong L(1,n^2)\oplus L(1,(n+1)^2).$ Similarly, let
$v^2$ be a highest weight vector of weight $(n+2)^2$ in $N.$ Then
$$v^2=u^2+aJ'_{2-2(n+1)}J'_{2-2n}v$$ for some $u^2\in N^0, 0\neq a\in\C$
and $N^1\cong L(1,n^2)\oplus L(1,(n+1)^2)\oplus L(1,(n+2)^2).$
Continuing in this way we show that $N^i\cong
\oplus_{j=0}^{i+1}L(1,(n+j)^2)$ for all $i\geq 1.$ This implies that
$N\cong {\bar N}$ as $L(1,0)$-modules. Consequently, $N=\bar{N}$ is
an irreducible $M(1)^+$-module. The proof is complete. \qed

For convenience, denote  $(M(1)^{+})^\bot_{k^2}$ by $B_{k}$. Let
$P^k$ be the $M(1)^{+}$-submodule of $V$ generated by $B_{k}$.

\begin{lem}\label{lem6.2}
 The restriction of $(\cdot,\cdot)$
to $P^k$ is still non-degenerate.
\end{lem}
\pf By (\ref{e6.1}), the restriction of $(\cdot,\cdot)$ to $B_{k}$
is non degenerate. By (\ref{e6.2}), $B_{k}$ is a
$\C[J'_{3}]$-module. Let $0\subset W^1\subset W^2 \subset
\cdots\subset W^t=B_k$ be a chain of $\C[J'_{3}]$-modules such that
$W^{i+1}/W^i$ is irreducible where $W^0=0.$ Let $S^i$ be the
$M(1)^+$-submodule of $V$ generated by $W^i.$ Then $S^t=P^k$ and
$S^{i+1}/S^i$ is an irreducible $M(1)^+$-module by Lemma
\ref{lem6.1}. This implies that if $u\in P^k$ satisfying
$L(m)u=J'_nu=0$ for all $m\geq 1$ and $n\geq 4$ and $L(0)u\in\C u,
J'_3u\in \C u$ then $u\in B_k.$

Now let $R$ be the radical of the restriction of the bilinear form
to $P^k.$ Then $R$ is an $M(1)^+$-submodule of $P^k.$ If $R\ne 0$,
then $R$ contains an irreducible $M(1)^+$-submodule whose
intersection with $B_k$ is nonzero. As a result the restriction of
the bilinear form to $B_k$ is degenerate, a contradiction. This
completes the proof. \qed

Using Lemma \ref{lem6.2} gives the following decomposition
$$
V=(M(1)^{+}\bigoplus P^k)\bigoplus (M(1)^{+}\bigoplus P^k)^\bot.
$$
Moreover, $(M(1)^{+}\bigoplus P^k)^\bot$ is an $M(1)^{+}$-module.
Then
$$
(M(1)^{+}\bigoplus P^k)^\bot=\bigoplus_{m\geq k_{1}}b_{m}L(1,m^2)
$$
as an $L(1,0)$-module where $k_{1}> k$ and $b_{k_{1}}\geq 1$.
Similar to Lemma \ref{lem6.2}, the restriction of $(\cdot,\cdot)$ to
the $M(1)^{+}$-module generated by $(M(1)^{+}\bigoplus
P^k)^\bot_{k_{1}^2}$ is non-degenerate. Continuing in this way we
deduce the following lemma.
\begin{lem}\label{lem6.3} As an $M(1)^{+}$-module, $V$ has the
following submodule decomposition:
$$V=M(1)^{+}\bigoplus(\bigoplus_{i=0}^{\infty}P^{k_{i}})$$ where
$k_{0}=k<k_{1}<k_{2}<\cdots $ and $P^{k_{i}}$ is the
$M(1)^{+}$-submodule of $V$ generated by
$$B_{k_{i}}=\{u\in
V_{k_{i}^2}|L_{m}u=J'_{n}u=0, m\geq 1, n\geq 4\}.$$
\end{lem}

We  can now prove the following lemma.
\begin{lem}\label{lem6.4} We have $\dim
B_{k_{i}}=1$ for $i=0,1,2,\cdots$.  Furthermore, $V$ is a completely
reducible $M(1)^{+}$-module.
\end{lem}
\pf By Lemma \ref{lem6.1} and Lemma \ref{lem6.3}, it is enough to
prove that $\dim B_{k_{i}}=1$, for $i\geq 0$. We only  prove $\dim
B_{k_0}=1$ as the proof for other cases is similar. Suppose $\dim
B_{k}\geq 2$. Then there exist $x^{1},x^2\in B_{k}$ such that
$(x^i,x^i)=0$ and $(x^1,x^2)=1$ for $i=1,2$ (see the proof of Lemma
5.2 of \cite{DJ2}),  and
\begin{equation}\label{e6.5}
J'_{3}x^1=(4k^4-k^2)x^1.
\end{equation}
Denote by $M^{i}$ the irreducible $L(1,0)$-submodule of $V$
generated by $x^i$ respectively, $i=1,2$. Then $M^1\cong M^2\cong
L(1,k^2)$.  We first prove

{\bf Claim:}  $M^1\cdot M^1\cong L(1,4k^2)$ is an irreducible
$L(1,0)$-module.

Let $N^1$ be the $M(1)^{+}$-submodule of $V$ generated by $x^1$.
Since $(x^1,x^1)=0$ we see that $M(1)^{+}\cap (N^1\cdot N^1)=0.$ By
Lemma \ref{lem6.3},
$$N^1\cdot N^1=\bigoplus_{i=0}^{\infty}(P^{k_{i}}\cap (N^1\cdot N^1))$$
and $P^{k_{i}}\cap (N^{1}\cdot N^1)$ is either zero or a direct sum
of indecomposable $M(1)^{+}$-modules with lowest weight $k_{i}^{2}$.
Then by Theorem \ref{2t2},
\begin{equation}\label{e6.6}P^{k_{i}}\cap (N^{1}\cdot N^1)=0,\ \ \ {\rm or} \
\  k_{i}^2=4k^2.
\end{equation}
This implies that $M^1\cdot M^1\subset N^1\cdot N^1\subset P^{2k}.$
Using the fusion rules from Theorem \ref{2t1} then forces $M^1\cdot
M^1\cong L(1,4k^2).$  So the claim holds.

The rest proof of the lemma is quite similar to that of Lemma 5.2 in
\cite{DJ2}. We omit it. \qed

We are now in a position to prove the main result of this paper.
\begin{theorem}
Let $V$ and $k$ be as above. Then $V$ is isomorphic to the rational
vertex operator algebra $V_{L}^{+}$, where $L=\Z\al$ is the rank one
positive definite even lattice  such that $(\al,\al)=2k^2$.

\end{theorem}
\pf By Lemma \ref{lem6.4}, $\dim B_{k}=1$, so there exists a
non-zero element $F^1\in V_{k^2}$ such that
$$
J'_{3}F^{1}=(4k^4-k^2)F^{1}, \ J'_{m}F^{1}=0, \ L(n)F^{1}=0, \ m\geq
4,n\geq 1,$$ and
$$(F^{1},F^{1})=2.$$
Let $V_{L}$ be the vertex operator algebra
 associated to the positive even lattice $L=\Z \al$ such that
 $(\al, \al)=2k$. Let $V^0$ be the vertex operator subalgebra of $V$ generated by
$F^{1},J',\omega$. We first prove that $V^0\cong V_{L}^{+}$.

For $m\in\Z_{+}$, set
$$E^{m}=e^{m\al}+e^{-m\al}\in V_{L}^{+}.$$
Then $$(E^m,E^m)=2.$$
 Denote by $N^m$ the $M(1)^{+}$-submodule of $V_{L}^{+}$
generated by $E^m$. Then (see \cite{DN1})
$$N^m=span\{u_{n}E^m|u\in M(1)^{+},n\in\Z\}$$
and
$$
N^m\cdot N^n=N^{m-n}\bigoplus N^{m+n}$$ for $m,n\in\Z_{+},m\geq n$
where $N^0=M(1)^{+}$.

Let $W^1$ be the $M(1)^{+}$-submodule of $V$ generated by $F^1$.
Then by Lemma \ref{lem6.1}, $W^{1}$ is irreducible. So as
$M(1)^{+}$-modules, $N^1\cong W^1$.  Let $\sigma$ be the
$M(1)^{+}$-module isomorphism from $N^1$ to $W^1$ such that
\begin{equation}\label{e6.7}
\sigma(E^1)=F^{1}.
\end{equation}
By Theorem \ref{2t2}, Lemma \ref{lem6.4} and the fact that
$(E^1,E^1)=(F^1,F^1)$, we have
$$
F^{1}_{n}F^{1}=E^1_{n}E^1\in M(1)^{+},
$$
for $n\geq -2k^2$. Note that
\begin{align*}
& (E^{1}_{-2k^2-1}E^1,E^{1}_{-2k^2-1}E^1)=(E^1,E^{1}_{4k^2-1}E^{1}_{-2k^2-1}E^1)\\
&
=(E^1,\sum\limits_{i=0}^{4k^2-1}\left(\begin{array}{c}4k^2-1\\i\end{array}\right)
(E^{1}_{j}E^1)_{2k^2-2-j}E^1),
\end{align*}
\begin{align*}
& (F^{1}_{-2k^2-1}F^1,F^{1}_{-2k^2-1}F^1)=(F^1,F^{1}_{4k^2-1}F^{1}_{-2k^2-1}F^1)\\
&
=(F^1,\sum\limits_{i=0}^{4k^2-1}\left(\begin{array}{c}4k^2-1\\i\end{array}\right)
(F^{1}_{j}F^1)_{2k^2-2-j}F^1).
\end{align*}
This implies that
$$
(E^{1}_{-2k^2-1}E^1,E^{1}_{-2k^2-1}E^1)=(F^{1}_{-2k^2-1}F^1,F^{1}_{-2k^2-1}F^1).$$

Assume that $E^1_{-2k^2-1}E^1=a_{1}E^2+u^1$, where $0\neq
a_{1}\in\C$, $u^1\in M(1)^{+}$. Then $F^1_{-2k^2-1}F^1-u^1$ is
either zero or a non-zero primary vector of weight $4k^2$. Since
$$(F^1_{-2k^2-1}F^1-u^1,F^1_{-2k^2-1}F^1-u^1)=(a_{1}E^1,a_{1}E^1),$$
it follows that $F^1_{-2k^2-1}F^1-u^1\neq 0$. Let $W^2$ be the
$M(1)^{+}$-submodule of $V$ generated by $F^1_{-2k^2-1}F^1-u^1$.
Then the isomorphism $\sigma$ can be extended to the isomorphism
from $N^1\bigoplus N^2$ to $W^1\bigoplus W^2$ such that
$$\sigma(E^1)=F^1,  \ \sigma(E^2)=F^2,$$
where $$F^2=\frac{1}{a_{1}}F^1_{-2k^2-1}F^1-u^1.$$ So for any
$n\in\Z,$
$$\sigma(E^1_{n}E^1)=(\sigma E^1)_{n}(\sigma E^1).$$

Following the proof of Lemma 5.7 of \cite{DJ2} and continuing in
this way we deduce that $V^0\cong V_{L}^{+}$. The rest of the proof
is the same as that of Theorem 5.8 of \cite{DJ2}. \qed

\section{Appendix}
\def\theequation{7.\arabic{equation}}
\setcounter{equation}{0}

Let $X^i$, $i=1,2,\cdots,s$ be the same as in Lemma \ref{l5.4}. Then
we have
\begin{lem}\label{app2}
\begin{align*}
&
 (X^{i}_{4}X^{j})_{2}X^{k}\\
= &
\dfrac{1}{2}a_{jk}^{i}L_{-1}X_{4}^{i}X^{i}+a_{ik}^{j}(\dfrac{28}{75}L_{-2}+\dfrac{11}{30}L_{-1}^{2})X^{j}_{5}X^{j}\\
&
-\dfrac{197}{150}a_{jk}^iL_{-1}X_{4}^{j}X^{j}+\dfrac{1}{27}(\dfrac{114}{75}a_{ik}^{j}-2a_{jk}^{i})X^{(0)}-
 2\sum\limits_{r=1}^{s}\sum\limits_{l=1}^{s}a_{jk}^{r}a_{ir}^{l}X^{l}\\
 & +\dfrac{114}{75}\sum\limits_{r=1}^{s}\sum\limits_{l=1}^{s}a_{ik}^{r}a_{jr}^{l}X^{l}+a\delta_{jk}X^{i}+b\delta_{ik}X^{j},
\end{align*}
for some $a,b\in\C$.
\end{lem}
\pf By the Jacobi identity, we have
\begin{align*}
&
 (X^{i}_{4}X^{j})_{2}X^{k}\\
= &
\sum\limits_{p=0}^4(-1)^{p}\left(\begin{array}{c}4\\p\end{array}\right)[X_{4-p}^iX_{2+p}^j-X^j_{6-p}X_{p}^i]X^k\\
= &
[X_{4}^iX_{2}^j-X_{6}^jX_{0}^i-4X_{3}^iX_{3}^j+4X_{5}^jX_{1}^i+6X_{2}^iX_{4}^j\\
&
-6X_{4}^jX_{2}^i-4X_{1}^iX_{5}^j+4X_{3}^jX_{3}^i-X_{2}^jX_{4}^i]X^k.
\end{align*}
 Since
$$
X_{3}^pX^{q}=\frac{1}{27}\delta_{pq}X^{(0)}+\sum\limits_{r=1}^{s}a_{pq}^{r}X^r,$$
we have
$$
X_{3}^iX_{3}^jX^k=X_{3}^i(\frac{1}{27}\delta_{jk}X^{(0)}+\sum\limits_{p=1}^sa_{jk}^{p}X^p)=
\frac{1}{27}\delta_{jk}X_{3}^iX^{(0)}+\frac{1}{27}a_{jk}^{i}X^{(0)}+\sum\limits_{p,q=1}^sa_{jk}^{p}a_{ip}^{q}X^q,
$$
$$
X_{3}^jX_{3}^iX^k=X_{3}^j(\frac{1}{27}\delta_{ik}X^{(0)}+\sum\limits_{p=1}^sa_{ik}^{p}X^p)=
\frac{1}{27}\delta_{ik}X_{3}^jX^{(0)}+\frac{1}{27}a_{ik}^{j}X^{(0)}+\sum\limits_{p,q=1}^sa_{ik}^{p}a_{jp}^{q}X^q.
$$
Recall Lemma \ref{app1} and Theorem \ref{2t1}. For $p,q\in
\{1,2,\cdots,s\}$,
$$
X_{2}^pX^q=\delta_{pq}u_{p,q}+\sum\limits_{r=1}^{s}a_{pq}^{r}\frac{1}{2}
L(-1)X^r,
$$
$$
X_{1}^pX^q=\delta_{pq}v_{p,q}+\sum\limits_{r=1}^{s}a_{pq}^{r}(\frac{28}{75}
L(-2)X^r+\frac{23}{300}L(-1)^{2}X^r),
$$
$$
X_{0}^pX^q=\delta_{pq}w_{p,q}+\sum\limits_{r=1}^{s}a_{pq}^{r}(\frac{14}{75}
L(-3)X^r+\frac{14}{75} L(-2)L(-1)X^r-\frac{1}{300} L(-1)^{3}X^r),$$
where $u_{p,q},v_{p,q},w_{p,q}\in L(1,0)$. Note that for $1\leq
p\leq s$,
$$
X^p_{4}L(-1)=L(-1)X^p_{4}+4X^p_{3},  \
X^p_{5}L_(-2)=L(-2)X^p_{5}+8X^p_{3},
$$
$$
X^p_{5}L(-1)^2=(L(-1)X^p_{5}+5X^p_{4})L(-1)=L(-1)^{2}X^p_{5}+10L(-1)X^p_{4}+20X^p_{3},$$$$
 X^p_{6}L(-3)=L(-3)X^p_{6}+12X^p_{3}, \
 X^p_{6}L(-1)=L(-1)X^p_{6}+6X^p_{5},$$
\begin{eqnarray*}
& &X^p_{6}L(-2)L(-1)=(L(-2)X^p_{6}+9X^p_{4})L(-1)\\
& &\ \ \ \ \
=L(-2)L(-1)X^p_{6}+6L(-2)X^p_{5}+9L(-1)X^p_{4}+36X^p_{3},
\end{eqnarray*}
$$
X^p_{6}L(-1)^{3}=
L(-1)^3X^p_{6}+18L(-1)^2X^p_{5}+90L(-1)X^p_{4}+120X^p_{3}.$$ Then we
have
\begin{align*}
& X_{4}^iX_{2}^jX^k\\
= &
\delta_{jk}X_{4}^iu_{j,k}+\sum\limits_{p=1}^{s}a_{jk}^{p}\frac{1}{2}X_{4}^iL(-1)X^p\\
= &
\delta_{jk}X_{4}^iu_{j,k}+\sum\limits_{p=1}^{s}a_{jk}^{p}\frac{1}{2}(L(-1)X_{4}^i+4X_{3}^i)X^p\\
=
&\delta_{jk}X_{4}^iu_{j,k}+\frac{1}{2}a_{jk}^{i}L(-1)X_{4}^iX^i+\frac{2}{27}a_{jk}^{i}X^{(0)}+
2\sum\limits_{p,q=1}^{s}a_{jk}^{p}a_{ip}^{q}X^q,
\end{align*}
\begin{align*}
& X_{4}^jX_{2}^iX^k\\
=
&\delta_{ik}X_{4}^ju_{i,k}+\frac{1}{2}a_{ik}^{j}L(-1)X_{4}^jX^j+\frac{2}{27}a_{ik}^{j}X^{(0)}+
2\sum\limits_{p,q=1}^{s}a_{ik}^{p}a_{jp}^{q}X^q,
\end{align*}
\begin{align*}
& X_{5}^jX_{1}^iX^k\\
= &
\delta_{ik}X_{5}^jv_{i,k}+X_{5}^j\sum\limits_{p=1}^{s}a_{ik}^{p}(\frac{28}{75}
L(-2)X^p+\frac{23}{300}L(-1)^{2}X^p)\\
=&
\delta_{ik}X_{5}^jv_{i,k}+\sum\limits_{p=1}^{s}a_{ik}^{p}\frac{28}{75}(L(-2)X_{5}^j+8X_{3}^j)X^p
\\& +\frac{23}{300}\sum\limits_{p=1}^{s}a_{ik}^{p}(L(-1)^{2}X_{5}^j+10L(-1)X_{4}^j+20X_{3}^j)X^p\\
= &
\delta_{ik}X_{5}^jv_{i,k}+\frac{28}{75}a_{ik}^{j}L(-2)X_{5}^jX^j+\frac{23}{300}a_{ik}^{j}L(-1)^{2}X_{5}^jX^j\\
& + \frac{23}{30}a_{ik}^{j}L(-1)X_{4}^jX^j+\frac{113}{25\times
27}a_{ik}^{j}X^{(0)}+\frac{113}{25}\sum\limits_{p,q=1}^{s}a_{ik}^{p}a_{jp}^{q}X^q,
\end{align*}
\begin{align*}
& X_{6}^jX_{0}^iX^k\\
= &
\delta_{ik}X_{6}^jw_{i,k}+X_{6}^j\sum\limits_{p=1}^{s}a_{ik}^{p}(\frac{14}{75}
L(-3)X^p+\frac{14}{75}L(-2)L(-1)X^p-\frac{1}{300}L(-1)^3X^p)\\
=&
\delta_{ik}X_{6}^jw_{i,k}+\sum\limits_{p=1}^{s}a_{ik}^{p}\frac{14}{75}(L(-3)X_{6}^j+12X_{3}^j)X^p
 \\ & +\frac{14}{75}\sum\limits_{p=1}^{s}a_{ik}^{p}(L(-2)L(-1)X_{6}^j+6L_{-2}X_{5}^j+9L(-1)X_{4}^j+36X_{3}^j)X^p\\
 &
 -\frac{1}{300}\sum\limits_{p=1}^{s}a_{ik}^{p}(L(-1)^3X_{6}^j+18L(-1)^2X_{5}^j+90L(-1)X_{4}^j+120X_{3}^j)X^p\\
= &
\delta_{ik}X_{6}^jw_{i,k}+\frac{28}{25}a_{ik}^{j}L(-2)X_{5}^jX^j-\frac{3}{50}a_{ik}^{j}L(-1)^{2}X_{5}^jX^j\\
& + \frac{69}{50}a_{ik}^{j}L(-1)X_{4}^jX^j+\frac{214}{25\times
27}a_{ik}^{j}X^{(0)}+\frac{214}{25}\sum\limits_{p,q=1}^{s}a_{ik}^{p}a_{jp}^{q}X^q.
\end{align*}
Note that
$$X_{2}^iX_{4}^jX^k,X_{1}^iX_{5}^jX^k, X^i_{4}\delta_{jk}u_{j,k}\in\C(\delta_{jk}X^i),$$
$$ X^j_{2}X_{4}^iX^k, X^j_{4}\delta_{ik}u_{i,k},
 X^j_{5}\delta_{ik}v_{i,k},X^j_{6}\delta_{ik}w_{i,k}
\in\C(\delta_{ik}X^j).$$  Then it is easy to check that the lemma
holds.\qed

\end{document}